\documentclass[12pt]{article} 
\def\Bbb{\bf} 
\newcommand\C{{ \Bbb C}}
\newcommand\R{{\Bbb R}}

\newcommand\g{{\gamma}}  
\newcommand\G{{\Gamma}} 
\newtheorem{lmm}{Lemma}
\newtheorem{thm}{Theorem}
\newtheorem{cor}{Corollary}

\newtheorem{prp}{Proposition}

\def\comment#1{ }

\newcommand{\dfrac}[2]{%
\frac{\displaystyle{#1}}{\displaystyle{#2}}}

\begin{document}
\title{
Intersection numbers of twisted cycles and the correlation functions of
the conformal field theory II\\}
\author{Katsuhisa MIMACHI \thanks{Department of Mathematics, Tokyo Institute of Technology, Tokyo 152-8551 Japan}\and 
\and Masaaki YOSHIDA \thanks{Department of Mathematics, Kyushu University, 
                             Fukuoka 810-8560 Japan}}
\maketitle
\noindent
{\textbf{Abstract:}}\ We evaluate the intersection numbers of loaded cycles and twisted forms associated with an $n$-fold Selberg-type integral. The result is deeply related with the geometry of the configuration space of $n+3$ points in the projective line. The polyhedral geometrical method developed in \cite{KY} and \cite{Yo2} is used in full. This is a continuation of \cite{MiY}.
\tableofcontents
\par\bigskip
\par\medskip\qquad\qquad\qquad\qquad\qquad\qquad{\bf Introduction}\par\bigskip
We evaluate the intersection numbers of loaded cycles associated with the $n$-fold integral -- a variant of the Selberg integral --
$$\int u(t)\ \frac{dt_1\wedge\cdots\wedge dt_n}{\prod_{i=1}^n t_i
(1-t_i) (z-t_i)},\quad u=\prod_{i=1}^n\ t_i^\alpha\
(1-t_i)^\beta\ (z-t_i)^\delta\ \prod_{1\le i<j\le n}(t_i-t_j)^{{2\gamma}},$$
where $z\not=0,1$ is a complex parameter.
Let $\mathcal S=\mathcal S_u$ be the local system defined by the integrand $u$ on
$$X=\{t=(t_1,\dots,t_n)\in\C^n\mid t_i\not=0,1,z,t_j\ (j\not=i)\ 1\le i\le n\},$$
$\mathcal S^-$ its dual, i.e. the local system on $X$ defined by $u^{-1}$.
When $z$ is real, the real locus $X_\R$ of $X$ breaks into disjoint $n$-cells. 
We load each cell in a standard way, say $\Delta$ with $u^{\pm1}$ as is explained in \S6,
 and make it loaded cycles $\Delta_\pm$, that is, elements of the {\it locally
finite} $n$-dimensional homology group $H_n^{\rm lf}(X,\mathcal S^\pm)$ with 
coefficients in $\mathcal S$. There is a natural dual pairing
$$H_n^{\rm lf}(X,\mathcal S)\times H_n(X,\mathcal S^-)\longrightarrow\C$$
called the {\it intersection pairing.} In this paper we assume that the 
exponents $\alpha,\beta,\gamma,\delta$ are sufficiently generic so that the natural map
$$H_n(X,\mathcal S)\longrightarrow H_n^{\rm lf}(X,\mathcal S)$$
is an isomorphism; the inverse map is called the {\it regularization} and is 
denoted by $\rm reg$.
We evaluate the intersection numbers of these cycles. 
The intersection number of the cycles $\Delta^1_+$ and ${\rm reg\ }\Delta^2_-$ 
will be denoted by $\Delta^1_+\bullet{\rm reg\ }\Delta^2_-$ or simply 
$\Delta^1\bullet\Delta^2$; $\Delta\bullet\Delta$ will often be called the 
self-intersection number of $\Delta$.
\par\medskip
The symmetric group $S_n$ acts on $X$ as permutations of the coordinates $t_j$.
Since the local system $\mathcal S$ is invariant under this action, 
$S_n$ acts also on the homology groups. We are specially interested in the 
symmetric part, which is by definition the subspace consisting of 
$S_n$-invariant elements. It is well? known that the symmetric part of 
$H_n^{\rm lf}(X,\mathcal S)$ is $(n+1)$-dimensional, 
and that (just for simplicity, we assume $1<z$)
$$\Delta^n_k:=\sum_{\sigma\in S_n}\{t\mid 0<t_{\sigma(1)}<\dots<t_{\sigma(k)}<1
, z<t_{\sigma(k+1)}<\cdots<t_{\sigma(n)}\},\quad k=0,\dots,n$$
form a basis. Since these are disjoint each other, among the 
intersection numbers $\Delta^n_i\bullet\Delta^n_j$, only the self-intersection
numbers only are non-zero. As is easily known from the expression above, 
$\Delta^n_k$ is the disjoint union of ${n\choose k}$ copies of the direct 
product $\Delta^k_k\times \Delta^{n-k}_{0}$. The self-intersection number 
of a direct product is the product of the self-intersections of the factors. 
Thus the computation in question boils down to that for $\Delta^n_0$ and 
$\Delta^n_n$, which are isomorphic. So in this paper we evaluate the self-intersection number of $\Delta^n=\Delta^n_n$.
\par\bigskip
In the complex $(x_1,\dots,x_n)$-space, we consider the multi-valued function
$$v:=\prod_{i=1}^n\ t_i^\alpha\ (1-t_i)^\beta\ \prod_{1\le i<j\le n}(t_i-t_j)^{{2\gamma}}.$$
Let $\mathcal S_v$ be the local system defined by $v$ on 
$$Y=\{(t_1,\dots,t_n)\in\C^n\mid t_i\not=0,1,t_j\ (j\not=i)\ 1\le i\le n\}.$$

\begin{thm}Let $D^\sigma$ be the loaded $n$-cycle loaded standardly with $v$
supported by
$$\{t=(t_1,\dots,t_n)\in Y_\R\mid 0<t_{\sigma(1)}<\dots<t_{\sigma(n)}<1\},\quad \sigma\in S_n$$
and put $\Delta^n:=\sum_{\sigma\in S_n}D^n_\sigma$. Then the self-intersection number $J_n:=\Delta^n\bullet \Delta^n$ is equal to
$$n!\prod_{j=1}^n\frac{1-abg^{n+j-2}}{(1-ag^{j-1})(1-bg^{j-1})}\ \frac{1-g}{1-g^j},$$
where $a=\exp2\pi i\alpha,\ b=\exp2\pi i\beta,\ g=\exp2\pi i\gamma$.
\end{thm}
\begin{thm} Let 
$$\omega:=\frac{dt_1\wedge\cdots\wedge dt_n}{\prod_{i=1}^n\ t_i\ (1-t_i)}$$
represent an element of $n$-th cohomology group $H^n(Y,\mathcal S_v)$  with coefficients in $\mathcal S_v$. 
The self-intersection number $\omega\bullet\omega$ is equal to 
$$(2\pi i)^n\prod_{j=1}^n\frac{\alpha+\beta+(n+j-2)\g}{(\alpha+(j-1)\g)(\beta+(j-1)\g)}.$$
\end{thm}

In this paper, {\it we give a direct geometric proof of Theorems 1 and 2.} 
In the complex $t=(t_1,\dots,t_n)$-space, we consider the hyperplanes
$$t_i=0,1,\ t_j\ (j\not=i),\quad i=1,\dots,n,$$
and blow-up the $t$-space along the non-normally crossing loci of the union of
these hyperplanes. Since the hyperplanes are defined over the reals, combinatorial information (e.g. intersection pattern) of the exceptional divisors together with these hyperplanes can be best understood if we study the real locus of the blew-up space: The real $t$-space is cut by these hyperplanes into $(n+2)!/2!$ (non-empty) pieces defined by
$$a_1<\cdots<a_{n+2},\quad \{a_1,\dots,a_{n+2}\}=\{0,t_1,\dots,t_n,1\},$$
among which are $n!$ bounded chambers (simplices) $\sigma(1)\cdots\sigma(n)$ defined by
$$0<t_{\sigma(1)}<\cdots<t_{\sigma(n)}<1,\quad \sigma\in S_n.$$
Blowing-up along the non-normal crossing loci in the complex space corresponds to the truncation of chambers in the real space. Each bounded chambers are truncated to be a $n$-polyhedron called `Terada-$n$'. This polyhedron is used firstly to describe homotopy $(n+2)$-associativity (in this context Terada-$n$ is called the Stasheff $n$-cell), and then to describe the $S_{n+3}$-equivariant minimal smooth compactification of the configuration space of the colored $n+3$ distinct points in the projective line (cf. \cite{Yo1}). Terada-1 is a segment, Terada-2 is a pentagon, Terada-3 is a polyhedron with six pentagonal faces and three rectangular faces, and Terada-$n$ is described in \S2.1. 

Let $T$ be the Terada-$n$ coded by $12\cdots n$ and $T^\sigma$ $(\sigma\in S_n$) the Terada-$n$ coded by $\sigma(1)\cdots\sigma(n)$. Once we know the adjacency of the $n!$ Terad-$n$'s, we can express the intersection number $T\bullet T^\sigma$ of $T$ and $T^\sigma$ loaded standardly with $v$ as a rational form in the exponents $a,b$ and $g$; it is the sum of the local contributions along {\it all faces} of $T\cap T^\sigma$. The sum
$$\frac{J_n}{n!}=\sum_{\sigma\in S_n}T\bullet T^\sigma$$
factorizes: The factors of the denominators correspond bijectively to the hyperfaces (with the exponents attached) of $T$, and the factors of the numerators correspond bijectively to the divisors (with the exponents attached) at infinity.

The self-intersection number of the $n$-form $\omega$ is the sum of the local contribution of {\it all vertices} of $T-\sum_{\sigma\in S_n}T^\sigma.$ Recall that evaluation of the intersection forms of cohomology is much simpler than that of homology, in general. Though it is also the case here, our geometric interpretation confirms that the cohomology intersection is the counterpart of the homology intersection.

\par\medskip\noindent
{\bf Relation with the Selberg integral.} The reader may notice the
resemblance between our Theorems  and the Selberg integral; indeed they are deeply related. Consider the Selberg function
$$\displaystyle{\mbox{\rm Sel}_n(\alpha,\beta,\g)}:=\displaystyle{\int_{(0,1)^n}\ t_i^{\alpha}\ (1-t_i)^{\beta}\
\prod_{1\le i<j\le n}|t_i-t_j|^{{2\gamma}}\ \prod_{i=1}^n \dfrac{dt_i}{t_i(1-t_i)}}.$$
If we admit the well-known fact that the symmetric part ($S_n$-invariant part) of $H_n(Y,\mathcal S_v)$ and $H^n(Y,\mathcal S_v)$ are 1-dimensional, the quadratic relation of the hypergeometric integral 
(see \cite{MaY}) implies the reciprocity relation
$$\mbox{Sel}_n(\alpha,\beta,\g)\cdot\mbox{Sel}_n(-\alpha,-\beta,-\g)=\Delta^n\bullet\Delta^n\ \cdot \ \omega\bullet\omega.$$
The Selberg formula
$$\mbox{\rm Sel}_n(\alpha,\beta,\g)=\displaystyle{\prod_{j=1}^n\frac{\G(\alpha+(j-1)\g)\ \G(\beta+(j-1)\g)\ \G(j\g+1)}{\G(\alpha+\beta+(n+j-2)\g)\ \G(\g+1)}}$$
and the reciprocity relation of the Gamma function
$$\G(\alpha)\G(1-\alpha)=\frac{\pi}{\sin \pi\alpha}\mbox{\quad or\quad}\G(\alpha)\G(-\alpha)=\frac{a}{1-a}\cdot\frac{2\pi i}{\alpha}\quad(a:=\exp2\pi i\alpha)$$
tell that Theorem 1 implies Teorem 2 and vice versa. The situation is exatly the same as the simplest case $n=1$: The reciprocity relation
$$B(\alpha,\beta)\cdot B(-\alpha,-\beta)=\frac{1-ab}{(1-a)(1-b)}\cdot\frac{2\pi i(\alpha+\beta)}{\alpha\beta}$$ 
of the Beta function 
$$B(\alpha,\beta):=\int_0^1t^\alpha(1-t)^\beta\frac{dt}{t(1-t)}$$
follows directly from the formula
$$B(\alpha,\beta)=\frac{\G(\alpha)\G(\beta)}{\G(\alpha+\beta)}$$
and the reciprocity relation of the Gamma function.

\section{Observation of 2D and 3D cases}
\subsection{2D: Geometry of a pentagon (\cite{MiY})}
In the $(x,y)$-space, we consider the lines
$$x=0,\ y=0,\quad x=1,\ y=1, \quad x=y,$$
with exponents $a,a,b,b,g^2$, respectively.
The non-normally crossing points are the two points
$$x=y=0,\quad x=y=1.$$
These lines cut the real $(x,y)$-plane into $4!/2!=12$ non-empty chambers 
$$t_0<t_1<t_2<t_3<t_4\quad \{t_0,t_1,t_2,t_3\}\subset\{0,x,y,1\},$$
among which are $2!$ bounded chambers defined by $0<x<y<1$ and $0<y<x<1$.
We blow-up the $(x,y)$-plane at the two singular points. Let us
consider the triangle $0xy1$ bounded by the three edges
$$0=x,\quad x=y, \quad y=1,$$
with the three vertices
$$0=x=y,\quad x=y=1,\quad \{0=x,\ y=1\}.$$
By the blowing-up, the triangle $0xy1$ is truncated to be a penatagon 
bounded by the five segments
$$0=x\ (a),\quad 0=x=y\ (a^2g^2),\quad  x=y\ (g^2),\quad  x=y=1\
(g^2b^2),\quad  y=1\ (b),$$
where the corresponding exponents are given in the parentheses. 
Here $0=x=y$ is regarded as the exceptional curve coming from the singular 
point $0=x=y$, through which three lines $0=x,0=y,x=y$ with exponents 
$a,a,g^2$ pass;
the exponent along this exceptional curve is just the product of these three 
exponents. We will also call this pentagon $0xy1$. 
As we explained in \cite{KY} the intersection number 
$0xy1\cdot 0xy1$ is the sum of 
the local contributions at the baricenters of all the possible
faces of the pentagonal cell $0xy1$, where the contribution of the
2-cell is 1, that of the 1-face $x=0$ is $\frac1{a-1}$, that of the 1-face
$y=1$ is $\frac1{c-1}$, that of the vertex $\{x=0\}\cap\{y=1\}$ is
$\frac1{(a-1)(c-1)}$, and so on. The intersection number $0xy1\cdot
0yx1$ of the two pentagons is $\frac{g}{g^2-1}$ times the intersection
number of the common edge. 
\begin{figure}
\begin{minipage}{6cm}
\unitlength 0.1in
\begin{picture}( 20.0000, 21.7000)( 21.5000,-23.6000)
%
\special{pn 8}%
\special{pa 2550 560}%
\special{pa 2550 2360}%
\special{fp}%
%
\special{pn 8}%
\special{pa 3550 2360}%
\special{pa 3550 560}%
\special{fp}%
%
\special{pn 8}%
\special{pa 3950 960}%
\special{pa 2150 960}%
\special{fp}%
%
\special{pn 8}%
\special{pa 2150 1960}%
\special{pa 3950 1960}%
\special{fp}%
%
\special{pn 8}%
\special{pa 3950 560}%
\special{pa 3950 560}%
\special{fp}%
\special{pa 3950 560}%
\special{pa 2150 2360}%
\special{fp}%
\put(23.5000,-15.6000){\makebox(0,0)[lb]{$a$}}%
\put(37.5000,-15.6000){\makebox(0,0)[lb]{$b$}}%
\put(29.5000,-21.6000){\makebox(0,0)[lb]{$a$}}%
\put(29.5000,-7.6000){\makebox(0,0)[lb]{$b$}}%
\put(41.5000,-3.6000){\makebox(0,0)[lb]{$g^2$}}%
\put(27.7000,-11.8000){\makebox(0,0)[lb]{$0xy1$}}%
\put(29.5000,-17.6000){\makebox(0,0)[lb]{$0yx1$}}%
\end{picture}%

\end{minipage}\qquad
\begin{minipage}{6cm}
\unitlength 0.1in
\begin{picture}( 26.0000, 22.0000)( 18.5000,-23.9000)
%
\special{pn 8}%
\special{pa 2050 1360}%
\special{pa 1850 1360}%
\special{fp}%
\special{sh 1}%
\special{pa 1850 1360}%
\special{pa 1918 1380}%
\special{pa 1904 1360}%
\special{pa 1918 1340}%
\special{pa 1850 1360}%
\special{fp}%
%
\special{pn 8}%
\special{pa 2450 560}%
\special{pa 2450 560}%
\special{fp}%
%
\special{pn 8}%
\special{pa 2450 560}%
\special{pa 2450 560}%
\special{fp}%
%
\special{pn 8}%
\special{pa 2450 560}%
\special{pa 2450 560}%
\special{fp}%
%
\special{pn 8}%
\special{pa 2450 560}%
\special{pa 3450 560}%
\special{fp}%
%
\special{pn 8}%
\special{pa 2650 360}%
\special{pa 2650 1360}%
\special{fp}%
\special{pa 3450 360}%
\special{pa 4450 1360}%
\special{fp}%
%
\special{pn 8}%
\special{pa 3850 560}%
\special{pa 3450 560}%
\special{fp}%
%
\special{pn 8}%
\special{pa 2650 1760}%
\special{pa 2650 1360}%
\special{fp}%
%
\special{pn 8}%
\special{pa 2450 1360}%
\special{pa 3450 2360}%
\special{fp}%
%
\special{pn 8}%
\special{pa 3050 2160}%
\special{pa 4450 2160}%
\special{fp}%
%
\special{pn 8}%
\special{pa 4250 960}%
\special{pa 4250 2360}%
\special{fp}%
%
\special{pn 8}%
\special{pa 4050 760}%
\special{pa 4050 760}%
\special{fp}%
%
\special{pn 8}%
\special{pa 2850 1960}%
\special{pa 4050 760}%
\special{fp}%
\put(30.5000,-3.6000){\makebox(0,0)[lb]{$b$}}%
\put(24.5000,-9.6000){\makebox(0,0)[lb]{$a$}}%
\put(34.5000,-25.6000){\makebox(0,0)[lb]{$a^2g^2$}}%
\put(42.5000,-5.6000){\makebox(0,0)[lb]{$g^2$}}%
\put(44.5000,-15.6000){\makebox(0,0)[lb]{$g^2b^2$}}%
\put(28.5000,-11.6000){\makebox(0,0)[lb]{$0xy1$}}%
\put(34.5000,-17.6000){\makebox(0,0)[lb]{$0yx1$}}%
%
\special{pn 8}%
\special{sh 1}%
\special{ar 2650 560 10 10 0  6.28318530717959E+0000}%
%
\special{pn 8}%
\special{sh 1}%
\special{ar 2650 560 10 10 0  6.28318530717959E+0000}%
\special{sh 1}%
\special{ar 2650 560 10 10 0  6.28318530717959E+0000}%
\put(28.5000,-21.6000){\makebox(0,0)[lb]{$a$}}%
\put(44.5000,-19.6000){\makebox(0,0)[lb]{$b$}}%
%
\special{pn 8}%
\special{ar 2660 570 0 0  0.0000000 6.2831853}%
%
\special{pn 8}%
\special{sh 1.000}%
\special{ar 2660 570 50 50  0.0000000 6.2831853}%
%
\special{pn 8}%
\special{sh 1.000}%
\special{ar 2650 1590 50 50  0.0000000 6.2831853}%
%
\special{pn 8}%
\special{sh 1.000}%
\special{ar 3630 560 50 50  0.0000000 6.2831853}%
\end{picture}%

\end{minipage}
\caption{Truncating triangles to be pentagons}
\end{figure}
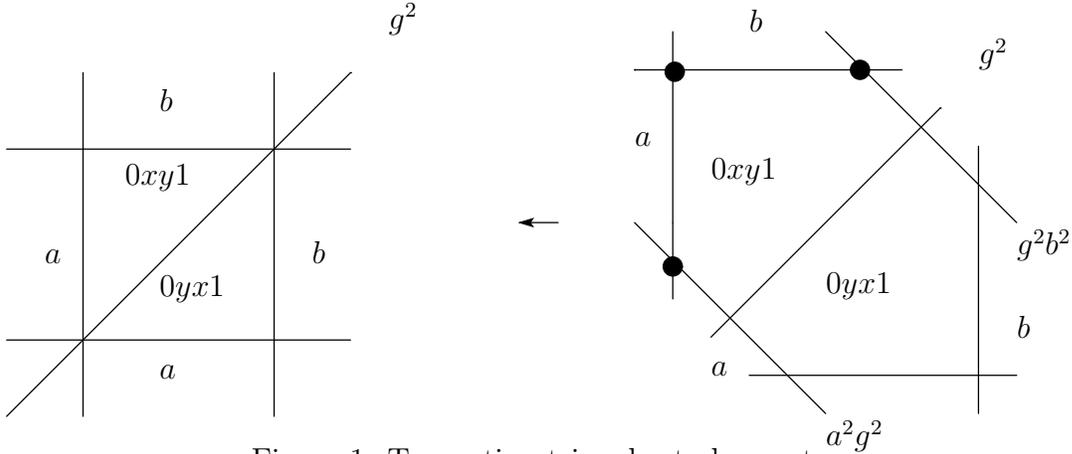
For notational simplicity, we introduce some functions:
$$F(x):=\frac1{x-1}, \quad S(x,y):=1+F(x)+F(y)\quad\left(=\frac{xy-1}{(x-1)(y-1)}\right),$$
$$\begin{array}{ll}
P(x_1,\dots,x_5):=
1&+F(x_1)+\cdots+F(x_5)\\
 &+F(x_1)F(x_2)+\cdots+F(x_4)F(x_5)+F(x_5)F(x_1).\end{array}$$
(Note that $S$ and $P$ are initials of segment and pentagon, while $F$ means nothing special.) Then we have
$$\begin{array}{ll}\displaystyle{\frac{J_2}2}&=0xy1\cdot(0yx1+0xy1)\\
&=-gF(g^2)S(a^2g^2,g^2b^2)+P(a,a^2g^2,a^2,g^2b^2,b).\end{array}$$
We evaluate this. Before directly computing
$$\begin{array}{ll}\displaystyle{\frac{J_2}2}&=\displaystyle{-\frac{g}{g^2-1}\left(\frac1{a^2g^2-1}+1+\frac1{g^2b^2-1}\right)}\\
&\displaystyle{+\ 1+\frac1{a-1}+\frac1{a^2g^2-1}+\frac1{g^2-1}+\frac1{g^2b^2-1}+\frac1{b-1}}\\
&\displaystyle{+\frac1{(a-1)(a^2g^2-1)}+\frac1{(a^2g^2-1)(g^2-1)}}\\
&\displaystyle{+\frac1{(g^2-1)(g^2b^2-1)}
+\frac1{(g^2b^2-1)(b-1)}+\frac1{(b-1)(a-1)}},
\end{array}$$
let us reflect a little:
The summands above tell that the donominator of $J_2$ divides
$$(a-1)(a^2g^2-1)(g^2-1)(b-1)(b^2g^2-1).$$
On the other hand, the exponents along the lines at infinity (after
blowing-up to make the singular lines cross normally) tell that the
numerator divides
$$(ag^2b-1)(a^2g^2b^2-1).$$
The numerator and the denominator have the same degree with respect to
each of $a,b$ and to $g$; they are symmetric with respect to $a$ and $b$. These
conditions are not enough to determine $J_2$; we do not know a priori
which factor is needed among $\{(g-1),(g+1)\}$, among
$\{(ag-1),(ag+1)\}$, among $\{(gb-1),(gb+1)\}$ and
among $\{(agb-1),(agb+1)\}$.

Let us add the terms above. First the terms free of $a$ and $b$:
$$1+\frac1{g^2-1}-\frac{g}{g^2-1}=\frac{g}{g+1}.$$
Next the terms with the factor $a^2g^2-1$ as denominators:
$$\frac1{a^2g^2-1}\left(1+\frac{1}{a-1}+\frac1{g^2-1}-\frac{g}{g^2-1}\right)
=\frac1{(a-1)(ag-1)(g+1)}.$$
Finally,
$$\begin{array}{ll}\displaystyle{\frac{J_2}2}&=\displaystyle{\frac{g}{g+1}+\frac1{a-1}+\frac1{b-1}+\frac1{(a-1)(b-1)}}\\
&\displaystyle{+\frac1{(a-1)(ag-1)(g+1)}+\frac1{(b-1)(gb-1)(g+1)}}\\
&\\
&\displaystyle{=\frac{abg+a+b-1}{(a-1)(b-1)(g+1)}+\frac1{(a-1)(ag-1)(g+1)}+\frac1{(b-1)(gb-1)(g+1)}};\end{array}$$
we know the denominator already. The numerator is a symmetric polynomial 
in $a$ and $b$; writing this as a polynomial in $ab$ and $a+b$, we can
easily check that the coefficients of the power of $a+b$ are $0.$
The result is
$$\frac{J_2}{2!}=\frac{(agb-1)(ag^2b-1)}{(a-1)(ag-1)(b-1)(gb-1)(g+1)}.$$
\subsection{3D: Geometry of Terada-3}
In the $(x,y,z)$-space, we consider the planes
$$x,\ y,\ z=0;\qquad x,\ y,\ z=1;\qquad x=y,\ y=z,\ z=x,$$
with exponents $a,b,g^2$, respectively.
The non-normally crossing loci are the lines
$$x=y=0,\ 1;\quad y=z=0,\ 1;\quad z=y=0,\ 1;\quad x=y=z$$
and the points
$$x=y=z=0,\ 1.$$
These planes cut the real $(x,y,z)$-space into $5!/2!$ chambers defined by
$$t_0<t_1<t_2<t_3<t_4\quad \{t_0,t_1,t_2,t_3,t_4\}\in\{0,x,y,z,1\},$$
among which are $3!$ bounded chambers defined by
$$0<t_1<t_2<t_3<1\quad \{t_1,t_2,t_3\}\in\{x,y,z\};$$
these six tetrahedra fill the cube $0<x,y,z<1$.
These chambers are encoded as $t_0\cdots t_4$; and $0t_1t_2t_31$ is often abbreviated as $t_1t_2t_3$.
Note that two such chambers are adjacent if and only if these codes can be changed from one to the other by a permutaion of adjacent two letters, 
and that two chambers $0xyz1$ and $01zyx$ are antipodal with center $(1,1,1)$.

We blow-up the $(x,y,z)$-space at the singular points and then along the singular lines. We describe the resulting object by the truncation of the chambers.
Let us consider the tetrahedron $0xyz1$ bounded by the four planes
$$0=x,\ x=y,\ y=z,\ z=1,$$
with the six edges 
$$0=x=y,\ x=y=z,\ y=z=1;\ \{0=x,y=z\},\ \{0=x,z=1\},\ \{x=y,z=1\},$$
and with the four vertices
$$0=x=y=z,\ x=y=z=1;\ \{0=x,y=z=1\},\ \{0=x=y,z=1\}.$$
The first three edges and the first two vertices above are in the singular loci. Such an edge is truncated to be a rectangle, and such a vertex is truncated to be a pentagon. As a result, the tetrahedron $xyz1$ is truncated to be a Terada-3, which is bounded by six pentagons and three rectangles.
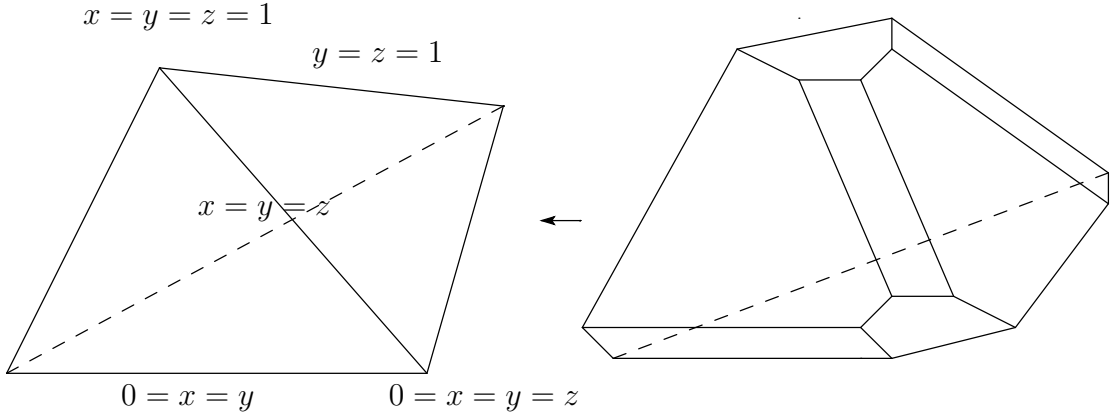
\begin{figure}
\begin{minipage}{6cm}
\unitlength 0.1in
\begin{picture}( 30.0000, 20.0000)( 17.8000,-24.3000)
%
\special{pn 8}%
\special{pa 2580 800}%
\special{pa 1780 2400}%
\special{fp}%
%
\special{pn 8}%
\special{pa 1780 2400}%
\special{pa 1780 2400}%
\special{fp}%
\special{pa 1780 2400}%
\special{pa 3980 2400}%
\special{fp}%
%
\special{pn 8}%
\special{pa 3980 2400}%
\special{pa 2580 800}%
\special{fp}%
%
\special{pn 8}%
\special{pa 2580 800}%
\special{pa 2580 800}%
\special{fp}%
\special{pa 2580 800}%
\special{pa 4380 1000}%
\special{fp}%
%
\special{pn 8}%
\special{pa 4380 1000}%
\special{pa 4380 1000}%
\special{fp}%
\special{pa 4380 1000}%
\special{pa 3980 2400}%
\special{fp}%
%
\special{pn 8}%
\special{pa 4380 1000}%
\special{pa 4380 1000}%
\special{da 0.070}%
\special{pa 4380 1000}%
\special{pa 1780 2400}%
\special{da 0.070}%
%
\special{pn 8}%
\special{pa 4780 1600}%
\special{pa 4780 1600}%
\special{fp}%
%
\special{pn 8}%
\special{pa 4780 1600}%
\special{pa 4780 1600}%
\special{fp}%
%
\special{pn 8}%
\special{pa 4780 1600}%
\special{pa 4780 1600}%
\special{fp}%
%
\special{pn 8}%
\special{pa 4780 1600}%
\special{pa 4780 1600}%
\special{fp}%
\special{pa 4780 1600}%
\special{pa 4580 1600}%
\special{fp}%
\special{sh 1}%
\special{pa 4580 1600}%
\special{pa 4648 1620}%
\special{pa 4634 1600}%
\special{pa 4648 1580}%
\special{pa 4580 1600}%
\special{fp}%
\put(33.8000,-8.0000){\makebox(0,0)[lb]{$y=z=1$}}%
\put(21.8000,-6.0000){\makebox(0,0)[lb]{$x=y=z=1$}}%
\put(23.8000,-26.0000){\makebox(0,0)[lb]{$0=x=y$}}%
\put(37.8000,-26.0000){\makebox(0,0)[lb]{$0=x=y=z$}}%
\put(27.8000,-16.0000){\makebox(0,0)[lb]{$x=y=z$}}%
\end{picture}%
\end{minipage}\qquad\qquad
\begin{minipage}{6cm}
\unitlength 0.1in
\begin{picture}( 27.5400, 17.8200)( 19.1000,-20.1200)
%
\special{pn 8}%
\special{pa 1910 1850}%
\special{pa 1910 1850}%
\special{fp}%
\special{pa 1910 1850}%
\special{pa 2072 2012}%
\special{fp}%
\special{pa 2072 2012}%
\special{pa 3530 2012}%
\special{fp}%
%
\special{pn 8}%
\special{pa 3530 2012}%
\special{pa 4178 1850}%
\special{fp}%
%
\special{pn 8}%
\special{pa 4178 1850}%
\special{pa 4664 1202}%
\special{fp}%
%
\special{pn 8}%
\special{pa 4664 1202}%
\special{pa 4664 1040}%
\special{fp}%
%
\special{pn 8}%
\special{pa 4664 1040}%
\special{pa 3530 230}%
\special{fp}%
%
\special{pn 8}%
\special{pa 1910 1850}%
\special{pa 3368 1850}%
\special{fp}%
%
\special{pn 8}%
\special{pa 3368 1850}%
\special{pa 3530 2012}%
\special{fp}%
%
\special{pn 8}%
\special{pa 3368 2012}%
\special{pa 3368 2012}%
\special{fp}%
\special{pa 3368 1850}%
\special{pa 3530 1688}%
\special{fp}%
%
\special{pn 8}%
\special{pa 3530 1688}%
\special{pa 3854 1688}%
\special{fp}%
%
\special{pn 8}%
\special{pa 3854 1688}%
\special{pa 4178 1850}%
\special{fp}%
%
\special{pn 8}%
\special{pa 3854 1688}%
\special{pa 3368 554}%
\special{fp}%
%
\special{pn 8}%
\special{pa 3368 554}%
\special{pa 3530 392}%
\special{fp}%
%
\special{pn 8}%
\special{pa 3530 392}%
\special{pa 3530 230}%
\special{fp}%
%
\special{pn 8}%
\special{pa 3530 392}%
\special{pa 4664 1202}%
\special{fp}%
%
\special{pn 8}%
\special{pa 3368 554}%
\special{pa 3044 554}%
\special{fp}%
%
\special{pn 8}%
\special{pa 3044 230}%
\special{pa 3044 230}%
\special{fp}%
\special{pa 3044 554}%
\special{pa 3530 1688}%
\special{fp}%
%
\special{pn 8}%
\special{pa 4664 1040}%
\special{pa 2072 2012}%
\special{da 0.070}%
%
\special{pn 8}%
\special{pa 2720 392}%
\special{pa 1910 1850}%
\special{fp}%
%
\special{pn 8}%
\special{pa 2720 392}%
\special{pa 3044 554}%
\special{fp}%
%
\special{pn 8}%
\special{pa 2720 392}%
\special{pa 3530 230}%
\special{fp}%
\end{picture}%

\end{minipage}
$$$$
\caption{Truncating a tetrahedron to be a Terada-3}
\end{figure}
This Terada-3 is adjacent to $0yxz1$ through the pentagon $x=y$, to $0zyx1$ through the (exceptional) rectangle born from the edge $x=y=z$, and to $01zyx$ through the (exceptional) pentagon born from the vertex $x=y=z=1$. These faces will be denoted by
$$\{xy\},\quad\{xyz\},\quad\{xyz1\},$$
respectively. In this way, we can easily know the codes of the adjacent (blown-up) chambers.

Once the faces of the truncated chambers and the adjacency through the faces are known, we can compute the intersection numbers following the recipe shown in 
\cite{KY}.

Let the exponents be given as follows:
$$x,y,z=0 \cdots a,\quad x,y,z=1 \cdots b,\quad x=y\cdots f^2,\quad y=z\cdots g^2,\quad z=x \cdots h^2.$$
The Terada-3 obtained from the tetrahedron $0xyz1$ 
has nine faces: the six penatagons
$$\{0x\},\ \{xy\},\ \{yz\},\ \{z1\},\ \{0xyz\},\ \{xyz1\},$$
and the three rectangles
$$ \{0xy\},\ \{xyz\},\ \{yz1\},$$
and 21 edges and 14 vertices (see Figure 3).
\begin{figure}\begin{minipage}{6cm}
\unitlength 0.1in
\begin{picture}( 31.8600, 29.8000)( 17.2800,-31.9000)
%
\special{pn 8}%
\special{pa 3352 788}%
\special{pa 2612 972}%
\special{fp}%
%
\special{pn 8}%
\special{pa 2612 972}%
\special{pa 2428 1712}%
\special{fp}%
\special{pa 2428 1712}%
\special{pa 2612 2452}%
\special{fp}%
%
\special{pn 8}%
\special{pa 4090 2452}%
\special{pa 4276 1712}%
\special{fp}%
%
\special{pn 8}%
\special{pa 4276 1712}%
\special{pa 4090 972}%
\special{fp}%
%
\special{pn 8}%
\special{pa 4090 972}%
\special{pa 3352 788}%
\special{fp}%
%
\special{pn 8}%
\special{pa 3352 788}%
\special{pa 3352 1158}%
\special{fp}%
%
\special{pn 8}%
\special{pa 2428 1712}%
\special{pa 2796 1712}%
\special{fp}%
%
\special{pn 8}%
\special{pa 3906 1712}%
\special{pa 4276 1712}%
\special{fp}%
%
\special{pn 8}%
\special{pa 4090 2452}%
\special{pa 3720 2266}%
\special{fp}%
%
\special{pn 8}%
\special{pa 2982 2266}%
\special{pa 2612 2452}%
\special{fp}%
%
\special{pn 8}%
\special{pa 2982 2266}%
\special{pa 2982 2266}%
\special{pa 2796 1712}%
\special{pa 3352 1158}%
\special{pa 3906 1712}%
\special{pa 3720 2266}%
\special{pa 2982 2266}%
\special{pa 2982 2266}%
\special{pa 2982 2266}%
\special{fp}%
\put(33.2000,-26.6600){\makebox(0,0){$\{0xyz\}$}}%
%
\special{pn 8}%
\special{pa 2612 2452}%
\special{pa 3352 2820}%
\special{fp}%
\special{pa 3352 2820}%
\special{pa 4090 2452}%
\special{fp}%
%
\special{pn 8}%
\special{pa 3352 2820}%
\special{pa 3352 3190}%
\special{fp}%
%
\special{pn 8}%
\special{pa 2612 972}%
\special{pa 2242 604}%
\special{fp}%
\put(33.2000,-19.8900){\makebox(0,0){$\{xy\}$}}%
\put(33.2000,-17.4200){\makebox(0,0){$g^2$}}%
\put(28.5800,-11.2600){\makebox(0,0){$b$}}%
\put(27.0400,-14.3400){\makebox(0,0){$\{z1\}$}}%
\put(37.8200,-11.2600){\makebox(0,0){$g^6b^3$}}%
\put(39.3600,-14.3400){\makebox(0,0){$\{xyz1\}$}}%
\put(33.2000,-3.5600){\makebox(0,0){$g^2b^2$}}%
\put(33.2000,-6.6400){\makebox(0,0){$\{yz1\}$}}%
\put(45.5200,-17.4200){\makebox(0,0){$g^2$}}%
\put(45.5200,-20.5000){\makebox(0,0){$\{yz\}$}}%
\put(20.8800,-17.4200){\makebox(0,0){$a$}}%
\put(20.8800,-20.5000){\makebox(0,0){$\{0x\}$}}%
\put(27.0400,-18.9600){\makebox(0,0){$a^2g^2$}}%
\put(39.3600,-22.0400){\makebox(0,0){$\{xyz\}$}}%
\put(40.9000,-18.9600){\makebox(0,0){$g^6$}}%
\put(33.2000,-23.5800){\makebox(0,0){$a^3g^6$}}%
%
\special{pn 4}%
\special{pa 3668 1474}%
\special{pa 2958 2182}%
\special{fp}%
\special{pa 3690 1496}%
\special{pa 2966 2220}%
\special{fp}%
\special{pa 3714 1520}%
\special{pa 2982 2252}%
\special{fp}%
\special{pa 3736 1542}%
\special{pa 3012 2266}%
\special{fp}%
\special{pa 3760 1566}%
\special{pa 3058 2266}%
\special{fp}%
\special{pa 3782 1588}%
\special{pa 3104 2266}%
\special{fp}%
\special{pa 3806 1612}%
\special{pa 3152 2266}%
\special{fp}%
\special{pa 3828 1636}%
\special{pa 3198 2266}%
\special{fp}%
\special{pa 3852 1658}%
\special{pa 3244 2266}%
\special{fp}%
\special{pa 3874 1682}%
\special{pa 3290 2266}%
\special{fp}%
\special{pa 3898 1704}%
\special{pa 3336 2266}%
\special{fp}%
\special{pa 3890 1758}%
\special{pa 3382 2266}%
\special{fp}%
\special{pa 3868 1828}%
\special{pa 3428 2266}%
\special{fp}%
\special{pa 3844 1896}%
\special{pa 3474 2266}%
\special{fp}%
\special{pa 3822 1966}%
\special{pa 3520 2266}%
\special{fp}%
\special{pa 3798 2036}%
\special{pa 3566 2266}%
\special{fp}%
\special{pa 3774 2104}%
\special{pa 3614 2266}%
\special{fp}%
\special{pa 3752 2174}%
\special{pa 3660 2266}%
\special{fp}%
\special{pa 3644 1450}%
\special{pa 2944 2150}%
\special{fp}%
\special{pa 3620 1428}%
\special{pa 2936 2112}%
\special{fp}%
\special{pa 3598 1404}%
\special{pa 2920 2082}%
\special{fp}%
\special{pa 3574 1380}%
\special{pa 2912 2044}%
\special{fp}%
\special{pa 3552 1358}%
\special{pa 2898 2012}%
\special{fp}%
\special{pa 3528 1334}%
\special{pa 2890 1974}%
\special{fp}%
\special{pa 3506 1312}%
\special{pa 2874 1944}%
\special{fp}%
\special{pa 3482 1288}%
\special{pa 2866 1904}%
\special{fp}%
\special{pa 3460 1266}%
\special{pa 2850 1874}%
\special{fp}%
\special{pa 3436 1242}%
\special{pa 2844 1836}%
\special{fp}%
\special{pa 3412 1220}%
\special{pa 2828 1804}%
\special{fp}%
\special{pa 3390 1196}%
\special{pa 2820 1766}%
\special{fp}%
%
\special{pn 4}%
\special{pa 3366 1174}%
\special{pa 2804 1736}%
\special{fp}%
%
\special{pn 4}%
\special{pa 4244 1820}%
\special{pa 3774 2290}%
\special{fp}%
\special{pa 4230 1882}%
\special{pa 3806 2304}%
\special{fp}%
\special{pa 4214 1944}%
\special{pa 3836 2320}%
\special{fp}%
\special{pa 4198 2004}%
\special{pa 3868 2336}%
\special{fp}%
\special{pa 4182 2066}%
\special{pa 3898 2352}%
\special{fp}%
\special{pa 4168 2128}%
\special{pa 3928 2366}%
\special{fp}%
\special{pa 4152 2190}%
\special{pa 3960 2382}%
\special{fp}%
\special{pa 4136 2252}%
\special{pa 3990 2398}%
\special{fp}%
\special{pa 4122 2312}%
\special{pa 4022 2412}%
\special{fp}%
\special{pa 4106 2374}%
\special{pa 4052 2428}%
\special{fp}%
\special{pa 4260 1758}%
\special{pa 3744 2274}%
\special{fp}%
\special{pa 4260 1712}%
\special{pa 3728 2244}%
\special{fp}%
\special{pa 4214 1712}%
\special{pa 3752 2174}%
\special{fp}%
\special{pa 4168 1712}%
\special{pa 3774 2104}%
\special{fp}%
\special{pa 4122 1712}%
\special{pa 3798 2036}%
\special{fp}%
\special{pa 4076 1712}%
\special{pa 3822 1966}%
\special{fp}%
\special{pa 4028 1712}%
\special{pa 3844 1896}%
\special{fp}%
\special{pa 3982 1712}%
\special{pa 3868 1828}%
\special{fp}%
\special{pa 3936 1712}%
\special{pa 3890 1758}%
\special{fp}%
%
\special{pn 8}%
\special{pa 4914 3168}%
\special{pa 4914 3168}%
\special{ip}%
%
\special{pn 8}%
\special{pa 4090 980}%
\special{pa 4860 210}%
\special{fp}%
%
\special{pn 8}%
\special{pa 4860 210}%
\special{pa 4860 3182}%
\special{ip}%
%
\special{pn 8}%
\special{pa 4860 3182}%
\special{pa 3166 3182}%
\special{ip}%
%
\special{pn 4}%
\special{pa 4860 2128}%
\special{pa 3806 3182}%
\special{fp}%
\special{pa 4860 2082}%
\special{pa 3760 3182}%
\special{fp}%
\special{pa 4860 2036}%
\special{pa 3714 3182}%
\special{fp}%
\special{pa 4860 1990}%
\special{pa 3668 3182}%
\special{fp}%
\special{pa 4860 1944}%
\special{pa 3620 3182}%
\special{fp}%
\special{pa 4860 1896}%
\special{pa 3574 3182}%
\special{fp}%
\special{pa 4860 1850}%
\special{pa 3528 3182}%
\special{fp}%
\special{pa 4860 1804}%
\special{pa 3482 3182}%
\special{fp}%
\special{pa 4860 1758}%
\special{pa 3436 3182}%
\special{fp}%
\special{pa 4860 1712}%
\special{pa 3390 3182}%
\special{fp}%
\special{pa 4060 2466}%
\special{pa 3352 3176}%
\special{fp}%
\special{pa 3968 2512}%
\special{pa 3352 3128}%
\special{fp}%
\special{pa 3874 2560}%
\special{pa 3352 3082}%
\special{fp}%
\special{pa 3782 2606}%
\special{pa 3352 3036}%
\special{fp}%
\special{pa 3690 2652}%
\special{pa 3352 2990}%
\special{fp}%
\special{pa 3598 2698}%
\special{pa 3352 2944}%
\special{fp}%
\special{pa 3506 2744}%
\special{pa 3352 2898}%
\special{fp}%
\special{pa 3412 2790}%
\special{pa 3352 2852}%
\special{fp}%
\special{pa 4860 1666}%
\special{pa 4098 2428}%
\special{fp}%
\special{pa 4860 1620}%
\special{pa 4114 2366}%
\special{fp}%
\special{pa 4860 1574}%
\special{pa 4130 2304}%
\special{fp}%
\special{pa 4860 1528}%
\special{pa 4144 2244}%
\special{fp}%
\special{pa 4860 1482}%
\special{pa 4160 2182}%
\special{fp}%
\special{pa 4860 1434}%
\special{pa 4176 2120}%
\special{fp}%
\special{pa 4860 1388}%
\special{pa 4190 2058}%
\special{fp}%
\special{pa 4860 1342}%
\special{pa 4206 1996}%
\special{fp}%
\special{pa 4860 1296}%
\special{pa 4222 1936}%
\special{fp}%
\special{pa 4860 1250}%
\special{pa 4236 1874}%
\special{fp}%
\special{pa 4860 1204}%
\special{pa 4252 1812}%
\special{fp}%
\special{pa 4860 1158}%
\special{pa 4268 1750}%
\special{fp}%
%
\special{pn 4}%
\special{pa 4860 1112}%
\special{pa 4276 1696}%
\special{fp}%
\special{pa 4860 1066}%
\special{pa 4260 1666}%
\special{fp}%
\special{pa 4860 1020}%
\special{pa 4252 1628}%
\special{fp}%
\special{pa 4860 972}%
\special{pa 4244 1588}%
\special{fp}%
\special{pa 4860 926}%
\special{pa 4236 1550}%
\special{fp}%
\special{pa 4860 880}%
\special{pa 4230 1512}%
\special{fp}%
\special{pa 4860 834}%
\special{pa 4214 1482}%
\special{fp}%
\special{pa 4860 788}%
\special{pa 4206 1442}%
\special{fp}%
\special{pa 4860 742}%
\special{pa 4198 1404}%
\special{fp}%
\special{pa 4860 696}%
\special{pa 4190 1366}%
\special{fp}%
\special{pa 4860 650}%
\special{pa 4182 1328}%
\special{fp}%
\special{pa 4860 604}%
\special{pa 4168 1296}%
\special{fp}%
\special{pa 4860 558}%
\special{pa 4160 1258}%
\special{fp}%
\special{pa 4860 510}%
\special{pa 4152 1220}%
\special{fp}%
\special{pa 4860 464}%
\special{pa 4144 1180}%
\special{fp}%
\special{pa 4860 418}%
\special{pa 4136 1142}%
\special{fp}%
\special{pa 4860 372}%
\special{pa 4122 1112}%
\special{fp}%
\special{pa 4860 326}%
\special{pa 4114 1072}%
\special{fp}%
\special{pa 4860 280}%
\special{pa 4106 1034}%
\special{fp}%
\special{pa 4860 234}%
\special{pa 4098 996}%
\special{fp}%
\special{pa 4860 2174}%
\special{pa 3852 3182}%
\special{fp}%
\special{pa 4860 2220}%
\special{pa 3898 3182}%
\special{fp}%
\special{pa 4860 2266}%
\special{pa 3944 3182}%
\special{fp}%
\special{pa 4860 2312}%
\special{pa 3990 3182}%
\special{fp}%
\special{pa 4860 2358}%
\special{pa 4036 3182}%
\special{fp}%
\special{pa 4860 2406}%
\special{pa 4082 3182}%
\special{fp}%
\special{pa 4860 2452}%
\special{pa 4130 3182}%
\special{fp}%
\special{pa 4860 2498}%
\special{pa 4176 3182}%
\special{fp}%
\special{pa 4860 2544}%
\special{pa 4222 3182}%
\special{fp}%
\special{pa 4860 2590}%
\special{pa 4268 3182}%
\special{fp}%
%
\special{pn 4}%
\special{pa 4860 2636}%
\special{pa 4314 3182}%
\special{fp}%
\special{pa 4860 2682}%
\special{pa 4360 3182}%
\special{fp}%
\special{pa 4860 2728}%
\special{pa 4406 3182}%
\special{fp}%
\special{pa 4860 2774}%
\special{pa 4452 3182}%
\special{fp}%
\special{pa 4860 2820}%
\special{pa 4498 3182}%
\special{fp}%
\special{pa 4860 2868}%
\special{pa 4544 3182}%
\special{fp}%
\special{pa 4860 2914}%
\special{pa 4592 3182}%
\special{fp}%
\special{pa 4860 2960}%
\special{pa 4638 3182}%
\special{fp}%
\special{pa 4860 3006}%
\special{pa 4684 3182}%
\special{fp}%
\special{pa 4860 3052}%
\special{pa 4730 3182}%
\special{fp}%
\special{pa 4860 3098}%
\special{pa 4776 3182}%
\special{fp}%
\special{pa 4860 3144}%
\special{pa 4822 3182}%
\special{fp}%
%
\special{pn 8}%
\special{sh 1.000}%
\special{ar 3344 804 38 38  0.0000000 6.2831853}%
%
\special{pn 8}%
\special{sh 1.000}%
\special{ar 2604 972 40 40  0.0000000 6.2831853}%
%
\special{pn 8}%
\special{sh 1.000}%
\special{ar 2442 1704 38 38  0.0000000 6.2831853}%
%
\special{pn 8}%
\special{sh 1.000}%
\special{ar 2620 2444 40 40  0.0000000 6.2831853}%
\put(27.1000,-21.9000){\makebox(0,0){$\{0xy\}$}}%
\end{picture}%

\end{minipage}\qquad\qquad\qquad\quad
\begin{minipage}{6cm}
The 2-faces of the Terada-3 $0xyz1$ are shown with their exponents.
The three 2-faces $$\{xy\},\ \{yz\},\ \{xyz\},$$ and the two edges
$$\{xy\}\cap\{xyz\},\ \{yz\}\cap\{xyz\}$$ touching five other Terada-3's are painted in gray, while the four vertices not touching any of these are marked.\end{minipage}
\caption{The faces of Terada-3 $0xyz1$ with different information}
\end{figure}

Let us evaluate the intersection numbers of the $S_3$ (acting as permutations of $\{x,y,z\}$) orbits of the simplex $0xyz1$.
The simplex $0xyz1$ (for simplicity we write $xyz$ for the time being) is truncated to be a Terada-3, which is also denoted by $xyz$. The Terada-3 touches other Terada-3s in the $S_3$-orbit along the following faces:
$$\begin{array}{lll}
xyz\cap yxz&=\mbox{pentagon}&\{xy\},\\
xyz\cap xzy&=\mbox{pentagon}&\{yz\},\\
xyz\cap zyx&=\mbox{rectangle}&\{xyz\},\\
xyz\cap zxy&=\mbox{segment}&\{xy\}\cap\{xyz\},\\
xyz\cap yzx&=\mbox{segment}&\{yz\}\cap\{xyz\}.\end{array}$$
We compute the intersection number
$$\frac{J_3}{3!}=xyz\cdot(yxz+xzy+zyx+zxy+yzx+xyz).$$
\begin{figure}
\begin{minipage}{6cm}
\unitlength 0.1in
\begin{picture}( 22.4000, 26.6200)( 20.8000,-28.4100)
%
\special{pn 8}%
\special{pa 3312 658}%
\special{pa 2640 826}%
\special{fp}%
%
\special{pn 8}%
\special{pa 2640 826}%
\special{pa 2472 1498}%
\special{fp}%
\special{pa 2472 1498}%
\special{pa 2640 2170}%
\special{fp}%
%
\special{pn 8}%
\special{pa 3984 2170}%
\special{pa 4152 1498}%
\special{fp}%
%
\special{pn 8}%
\special{pa 4152 1498}%
\special{pa 3984 826}%
\special{fp}%
%
\special{pn 8}%
\special{pa 3984 826}%
\special{pa 3312 658}%
\special{fp}%
%
\special{pn 8}%
\special{pa 3312 658}%
\special{pa 3312 994}%
\special{fp}%
%
\special{pn 8}%
\special{pa 2472 1498}%
\special{pa 2808 1498}%
\special{fp}%
%
\special{pn 8}%
\special{pa 3816 1498}%
\special{pa 4152 1498}%
\special{fp}%
%
\special{pn 8}%
\special{pa 3984 2170}%
\special{pa 3648 2002}%
\special{fp}%
%
\special{pn 8}%
\special{pa 2976 2002}%
\special{pa 2640 2170}%
\special{fp}%
%
\special{pn 8}%
\special{pa 2976 2002}%
\special{pa 2976 2002}%
\special{pa 2808 1498}%
\special{pa 3312 994}%
\special{pa 3816 1498}%
\special{pa 3648 2002}%
\special{pa 2976 2002}%
\special{pa 2976 2002}%
\special{pa 2976 2002}%
\special{fp}%
%
\special{pn 8}%
\special{pa 2640 2170}%
\special{pa 3312 2506}%
\special{fp}%
\special{pa 3312 2506}%
\special{pa 3984 2170}%
\special{fp}%
%
\special{pn 8}%
\special{pa 3312 2506}%
\special{pa 3312 2842}%
\special{fp}%
%
\special{pn 8}%
\special{pa 2640 826}%
\special{pa 2304 490}%
\special{fp}%
%
\special{pn 8}%
\special{pa 3984 826}%
\special{pa 4320 490}%
\special{fp}%
\put(31.4400,-3.4900){\makebox(0,0)[lb]{$r_1$}}%
\put(27.2400,-10.7000){\makebox(0,0)[lb]{$p_2$}}%
\put(35.5000,-10.5600){\makebox(0,0)[lb]{$p_1$}}%
\put(31.7900,-15.9500){\makebox(0,0)[lb]{$p_3$}}%
\put(25.8400,-18.5400){\makebox(0,0)[lb]{$r_2$}}%
\put(37.6000,-18.7500){\makebox(0,0)[lb]{$r_3$}}%
\put(43.1300,-22.8800){\makebox(0,0)[lb]{$q_1$}}%
\put(32.0000,-22.9500){\makebox(0,0)[lb]{$q_2$}}%
\put(20.8000,-22.9500){\makebox(0,0)[lb]{$q_2$}}%
\end{picture}%

\end{minipage}
\qquad\qquad
\begin{minipage}{6cm}
The nine faces of a Terada-3 are shown with exponts $p_i,$ $q_i,$
 $r_i$. The self-intersection number of this Terada-3 is given by $$Q(p_1,q_1,r_1;p_2,q_2,r_2;p_3,q_3,r_3).$$
\end{minipage}
\caption{A stereographic image of a Terada-3}
\end{figure}
Notation for Terada-3: $Q(p_1,q_1,r_1;p_2,q_2,r_2;p_3,q_3,r_3)$
$$=Q\left(\begin{array}{ccc}p_1&&q_1\\&r_1&\\p_2&&q_2\\&r_2&\\p_3&&q_3\\&r_3&\end{array}\right)=\begin{array}{l}1+\sum F(r_i)S(p_i,q_{i+1})S(p_{i+1},q_i)\\\\+\sum F(p_i)+F(q_i)+ F(p_i)F(p_{i+1})+F(q_i)F(q_{i+1})\\\\+F(p_1)F(p_2)F(p_3)+F(q_1)F(q_2)F(q_3)+\sum F(p_i)F(q_i).\end{array}$$
(summation on $i$ is modulo 3.) Then, taking account of orientation (see \cite{KY} \S3), we have 
$$\begin{array}{ll}
J_3/3!&=fF(f^2)P(b,a^2f^2,a^3f^2g^2h^2,f^2g^2h^2,f^2g^2h^2b^3)
\quad\mbox{pentagon $yxz$}\\[1mm]
&+gF(g^2)P(a,g^2b^2,f^2g^2h^2b^3,f^2g^2h^2,a^3f^2g^2h^2) 
\quad\mbox{pentagon $xzy$}\\[1mm]
&-fghF(f^2g^2h^2)S(a^3f^2g^2h^2,f^2g^2h^2b^3)S(f^2,g^2) 
\quad\mbox{rectangle $zyx$}\\[1mm]
&+fghF(f^2g^2h^2)fF(f^2)S(a^3f^2g^2h^2,f^2g^2h^2b^3) 
\quad\mbox{segment $zxy$}\\[1mm]
&+fghF(f^2g^2h^2)gF(g^2)S(a^3f^2g^2h^2,f^2g^2h^2b^3) 
\quad\mbox{segment $yzx$}\\[1mm]
&-Q(f^2g^2h^2b^3, g^2 g^2b^2;b, a,a^2f^2; f^2, a^3f^2g^2h^2,f^2g^2h^2)
, \quad\mbox{$xyz$ itself}\end{array}$$
putting $f=h=g$, this is equal to
$$\begin{array}{ll}
&gF(g^2)P(b,a^2g^2,a^3g^6,g^6,g^6b^3)
+gF(g^2)P(a,g^2b^2,g^6b^3,g^6,a^3g^6) 
\quad\mbox{pentagons}\\[1mm]
&-g^3F(g^6)S(a^3g^6,g^6b^3)S(g^2,g^2) 
\quad\mbox{rectangle}\\[1mm]
&+g^3F(g^6)gF(g^2)S(a^3g^6,g^6b^3) 
+g^3F(g^6)gF(g^2)S(a^3g^6,g^6b^3) 
\quad\mbox{segments}\\[1mm]
&-Q(g^6b^3,g^2,g^2b^2;b,a,a^2g^2;g^2,a^3g^6,g^6),\quad\mbox{$xyz$ itself}
\end{array}$$
where the last term $Q$ is equal to $1+$
$$\begin{array}{ll}
&F(g^2b^2)S(g^2,b)S(a,g^6b^3)+F(a^2g^2)S(g^2,a)S(b,a^3g^6)
+F(g^6)S(g^2,g^2)S(a^3g^6,g^6b^3)\\[1mm]
&+F(g^2)+F(g^2)+F(a)+F(b)+F(a^3g^6)+F(g^6b^3)\\[1mm]
&+F(a^3g^6)(F(g^2)+F(a)+F(g^2))
+F(g^6b^3)(F(g^2)+F(b)+F(g^2))\\[1mm]
&+F(a)F(g^2)+F(b)F(g^2)+F(a)F(b)
+F(a)F(g^2)F(a^3g^6)+F(b)F(g^2)F(g^6b^3).\end{array}$$
Let us add these terms. First the terms free of $a$ and $b$: 
$$\begin{array}{ll}
2gF(g^2)-1-2F(g^2)\\[1mm]
+F(g^6)[\ 2gF(g^2)+2g^4F(g^2)-g^3S(g^2,g^2)-S(g^2,g^2)\ ]\\[1mm]
=\displaystyle{1-\frac{2g}{g+1}-\frac1{(g+1)(g^2+g+1)}}.\end{array}$$
Next the terms with the factor $a^3g^6-1$ as denominators: $(a^3g^6-1)^{-1}$ times
$$\begin{array}{ll}F(a^2g^2)[\ gF(g^2)-S(g^2,a)\ ]+F(g^6)[\ \mbox{as above}\ ]\\[1mm]
+F(a)(gF(g^2)-F(g^2)-1)+2gF(g^2)-2F(g^2)-1\\[1mm]
=\displaystyle{\frac{a^2g^4+ag^2+1}{(a-1)(ag-1)(g+1)(g^2+g+1)}}.\end{array}$$
Among the remaining terms, the terms with the factor $a^2g^2-1$ as denominators: $(a^2g^2-1)^{-1}$ times
$$\left(1+\frac1{b-1}\right)\frac{ag+1}{(a-1)(g+1)}.$$
Finally left are
$$F(a)(gF(g^2)-F(g^2)-1)+F(b)(gF(g^2)-F(g^2)-1)-F(b)F(a).$$
So we know the denominator of $J_3$; the numerator can be computed without difficulty and we get 
$$J_3=-\frac{3!\ (g^2ba - 1) (g^3ba - 1)(g^4b a - 1) }
{(a - 1) (ag - 1)(ag^2 - 1) (b - 1)(g b - 1)  (g^2b - 1)(g + 1) (g^2+ g + 1) }$$

\section{$n$D case}
In the real $(x_1,\dots,x_n)$-space, we consider the hyperplanes
$$x_i=0,1,\ x_j\ (j\not=i),\quad i=1,\dots,n.$$
They cut the space into $(n+2)!/2!$ chambers (simplices) defined by
$$a_1<\cdots<a_{n+2},\quad \{a_1,\dots,a_{n+2}\}=\{0,x_1,\dots,x_n,1\},$$
among which are $n!$ bounded chambers defined by
$$0<x_{\sigma(1)}<\cdots<x_{\sigma(n)}<1,\quad \sigma\in S_n.$$
Now we blow-up along the non-normally crossing loci of the union of
these hyperplanes. Each simplex is
transformed into an $n$-polygon called Terada-$n$. (A Terada-2
is a pentagon, and a Terada-3 appeared in the previous subsection.)  We
encode the Terada-$n$ coming from the simplex
$0<x_{\sigma(1)}<\cdots<x_{\sigma(n)}<1$ by the word 
$$0\sigma(1)\cdots\sigma(n)(n+1)(n+2),$$
which will be called an $(n+3)$-{\it juzu};
note that point $1$ is coded by $n+1$. 
\subsection{Shape of Terada-$n$ and Juzus}
A $k$-juzu is a sequence of $k$ numerals coding a Terada-$(k-3)$; two
juzus are identified if one is obtained from the other by a cyclic
permutation and/or inversion. For example, a $3$-juzu
$$012=120=\cdots=210$$ 
is a point, a $4$-juzu
$$0123=1230=\cdots=3210$$ 
is a segment, and a $5$-juzu
is a pentagon. To make the formulae looks nicer, we
often abbreviate $n\pm1$ and $n\pm2$ as
$$n-2={}''n,\quad n-1={}'n,\quad n+1=n',\quad n+2=n''.$$ Let us describe the boundary of the Terada-$n$ $T:=01\cdots\ nn'n''$
coming from the $n$-simplex
$$0<x_1<\cdots<x_n<1.$$ The Terada-$(n-k)$ coming from the
$(n-k)$-simplex defined by
$$0<x_1<\cdots<x_{i-1}<x_i=\cdots=x_{i+k}<x_{i+k+1}<\cdots<x_n<1,$$
which is part of the boundary of the $n$-simplex above, will be coded by
the $(n-k+3)$-juzu
$$0\ 1\cdots {}'i\ \left(i\ \cdots\phantom{{}^R} i+k\right)\ (i+k)'\ \cdots\ n\ n'\ n'',$$
regarding the $k+1$ numerals in the parentheses as a single numeral.
The $(n-1)$-faces of the Terada-$n$ $T=01\cdots\ nn'n''$ are direct
products of Terada-$(n-i)$ and Terada-$(i-1)$:
$$\begin{array}{l}
(0\cdots i)\ i'\cdots n''\times0\cdots i(i'\cdots n''),\quad
(1\cdots i')\ i''\cdots n''0\times1\cdots i'(i''\cdots n''0),\dots\\[1mm]
(n''0\cdots {}'i)\ i\cdots n'\times n''0\cdots {}'i\ (i\cdots n').\quad i=1,\dots,\left[\dfrac{n+3}2\right];
\end{array} $$
they will be abbreviated as
$$(0\cdots i)\mbox{ or }(i'\cdots n''),\quad (1\cdots i')\mbox{ or
}(i''\cdots n''0),\ \dots,\ 
(n''0\cdots {}'i)\mbox{ or }(i\cdots n'), $$
respectively.
The $(n-k)$-faces are intersections of $k$ of these $(n-1)$-faces, say
$$(a_1a'_1\cdots z_1),\ (a_2a'_2\cdots z_2),\ \dots,\ (a_ka'_k\cdots z_k)$$ 
such that for any $i$ and $j$, $\{a_i,a'_i,\dots,z_i\}\not=\{a_j,a'_j,\dots,z_j\}$
and either
$$\{a_i,a'_i,\dots,z_i\}\cap\{a_j,a'_j,\dots,z_j\}=\emptyset\ \ \mbox{ or}$$
$$ \{a_i,a'_i,\dots,z_i\}\subset\{a_j,a'_j,\dots,z_j\}\ \mbox{ or
}\ \{a_j,a'_j,\dots,z_j\}\subset\{a_i,a'_i,\dots,z_i\}.$$ 
{\bf Terminology:} In the sequel, in place of saying `for any $i$ and $j$, $\dots$', we say `the sets $$\{a_1,a'_1,\dots,z_1\},\ \{a_2,a'_2,\dots,z_2\},\ \dots,\ \{a_k,a'_k,\dots,z_k\}$$ have the {\it disjoint/include property}.'\par\smallskip
Note that this
definition fits the two alternative expression of an $(n-1)$-face above.
\par\noindent
\begin{itemize}
\item When $n=1$, Terada-1 $0123$ is a segment. Its boundary consists of two
points: $$(01)=(23),\quad (12)=(30).$$
\item When $n=2$, Terada-2 $012345$ is a pentagon. Its boundary consists of
\begin{itemize}
\item five edges: $$(01),\  (12),\  (34),\  (45),\  (50)$$ 
\item and five vertices
$$(01)\cap(34),\ (34)\cap(51),\ (50)\cap(12),\  (12)\cap(45),\  (45)\cap(01).$$
\end{itemize}
\item When $n=3$, the boundary of the Terada-3 $012345$ consists of 
\begin{itemize}
\item 2D faces of 2 kinds
$$(01),\ \dots,\ (50);\quad (012)=(345),\ (123)=(456),\ (234)=(501),$$
\item 1D faces of 2 kinds
$$(01)\cap(23);\quad (01)\cap(012),$$
\item 0D faces of 2 kinds
$$(01)\cap(23)\cap(45),\ (12)\cap(34)\cap(50);\quad (01)\cap(12)\cap(012).$$
\end{itemize}\end{itemize}
Let us see how this Terada-$n$ $T$ touches other Terada-$n$'s.
Through every face
$$(a_1a'_1\cdots z_1)\cap(a_2a'_2\cdots z_2)\cap\cdots\cap(a_ka'_k\cdots z_k)$$
with
$$a_1,a'_1,\dots,z_1;\ a_2,a'_2,\dots,z_2;\quad\dots\quad; a_k,a'_k,\dots,z_k\in\{1,\dots,n\}$$ 
touches $T$ the Terada-$n$ 
$$0\sigma(1)\cdots\sigma(n)n'n'',$$
where $\sigma\in S_n$ is the product of the {\it commutative} $k$ cyclic
permutaion of elements 
$$\{a_1,a'_1,\dots,z_1\},\ \{a_2,a'_2,\dots,z_2\},\ \dots,\ \{a_k,a'_k,\dots,z_k\}.$$
Note that not all of the $n!$ Terada-$n$'s touches the Terada-$n$ $T$.
\begin{itemize}
\item When $n=2$, the Terada-2 $01234$ touches the Terada-2 $02134$ through the
edge $(12)$. Note that if you perform the cyclic permutaion of
$\{3,4,0\}$ you get also the juzu $02134$.\par\noindent
\item When $n=3$, the Terada-3 $012345$ touches the five others 
\begin{itemize}
\item through 2D faces: $(12),\  (23);\quad (123),$
\item through 1D faces: $(12)\cap(123),\ (23)\cap(123).$
\end{itemize}
\item When $n=4$ already, among the $4!$ Terada-4's, there are Terada-4's which
do not touch $0123456$; they are $0241356$ and $0314256$.
\end{itemize}
\subsection{Exponents along the hyperfaces}
Since 
$$u=\prod_{i=1}^n\ t_i^\alpha\ (1-t_i)^\beta\ \prod_{1\le i<j\le n}(t_i-t_j)^{2\gamma},$$
the exponents along the hyperplanes in the $t$-space are given as 
follows
$$t_i=0\ \cdots\ a,\quad t_i=1\ \cdots\ b, \quad x_i=x_j\ \cdots\ g^2,$$
where $a=\exp2\pi i\alpha,\ b=\exp2\pi i\beta,\ g=\exp2\pi i\gamma$.
Through the $(n-k)$-simplex 
$$0=x_1=\cdots=x_{k}<x_{k+1}<\cdots<x_n<1,$$
pass the hyperplanes
$$0=x_i,\quad 1\le i\le k \mbox{\quad and \quad} x_i=x_j,\quad 1\le i<j\le k;$$
through the $(n-k)$-simplex 
$$0<x_1<\cdots<x_{p-1}<x_p=\cdots=x_{p+k}<x_{p+k+1}<\cdots<x_n<1,$$
pass the hyperplanes
$$x_i=x_j,\quad p\le i<j\le p+k;$$
through the $(n-k)$-simplex 
$$0<x_1<\cdots<x_{n-k}<x_{n-k+1}=\cdots=x_n=1,$$
pass the hyperplanes
$$ x_i=x_j,\quad n-k+1\le i<j\le n\mbox{\quad and \quad} x_i=1,
\quad n-k+1\le i\le n.$$
Thus, after necessary blowing-up, the exponents along the hyperfaces
 ($(n-1)$-faces)  of the Terada-$n$ $T=01\dots nn'n''$ are give as
$$(01\cdots k)\ \cdots\ a^kg^{k(k-1)},\quad (p\cdots p+k)\ \cdots\ g^{k(k+1)},
\quad (n-k+1\cdots n')\ \cdots\ b^kg^{k(k-1)},$$
where $1\le k,p,k+p\le n.$ 
\subsection{Intersection numbers}
The self-intersection number of an $n$-polygon bounded by {\it normally crossing} hyperplanes is the sum of 
the local contributions of all the possible
faces, where the contribution of the $n$-face (the plygon itself) is $1$, that of a hyperface with exponent $e$ is 
$\frac1{e-1}$, and that of a $p$-codimensional face (which is the intersection of exactly $p$ hyperfaces) is the product of the contributions of the $p$ intersecting hyperfaces.
\par\medskip\noindent
{\bf Notation:} For a hyperface $(qq'\cdots r)\ 0\le q<r\le n'$ with exponent $e$, we set
$$[q\cdots r]:=\dfrac1{e-1}.$$
Thus, taking the Terada-$n$ as a plygon, we have
\begin{prp}The self-intersection number $T\bullet T$ of the 
{\rm Terada-}$n$ $T$ is the sum of all the possible product
$$[a_1a'_1\cdots z_1][a_2a'_2\cdots z_2]\cdots[a_ka'_k\cdots z_k],$$ 
$0\le a_1,\dots,z_1,a_2,\dots,z_2,\cdots,a_k,\dots,z_k\le n'=n+1,$ $0\le k\le n'$, 
such that $[a_ia'_i\cdots,z_i]\not=[01\cdots n']$ and,
 $$\{a_i,a'_i,\dots,z_i\}\quad i=1,\dots,k$$ have the disjoint/include property.  When $k=0$ the 
product is regarded as $1$.
\end{prp}
\begin{prp}Let a {\rm Terada-}$n$ $T'$ touches $T$ along 
the hyperface $H:=(p,p',\dots,p+q)$ $(1\le p,p+q\le n)$. 
The intersection number of the Terada-$n$'s $T$ and $T'$ is given by
$$T\bullet T'=(-)^qg^{q(q+1)/2}[pp'\cdots p+q]\ H\bullet H,$$
where the self-intersection number $H\bullet H$ of $H$ is the sum of 
all the possible product
$$[a_1a'_1\cdots z_1][a_2a'_2\cdots z_2]\cdots[a_ka'_k\cdots z_k],$$ 
$0\le a_1,\dots,z_1,a_2,\dots,z_2,\dots,a_k,\dots,z_k\le n'=n+1,$ $0\le k\le n'$, 
such that $[a_ia'_i\cdots z_i]\not=[01\cdots n']$ and,
$$\{p,p',\dots,p+q\},\quad \{a_i,a'_i,\dots,z_i\}\quad i=1,\dots,k$$ have the disjoint/include property.
\end{prp}
\begin{cor}Let a {\rm Terada-}$n$ $T'$ touches $T$ along 
a face $F$ which is the intersection of the hyperfaces 
$H_l:=(p_l,p'_l,\dots,p_l+q_l)$ $(1\le p_l,p_l+q_l\le n)$. 
The intersection number of the {\rm Terada-}$n$'s $T$ and $T'$ is given by
$$T\bullet T'=\prod_l(-)^{q_l}g^{q_l(q_l+1)/2}[p_lp'_l\cdots p_l+q_l]\ F\bullet F,$$
where the self-intersection number $F\bullet F$ of $F$ is the sum of 
all the possible product
$$[a_1a'_1\cdots z_1][a_2a'_2\cdots z_2]\cdots[a_ka'_k\cdots z_k],$$ 
$0\le a_1,\dots,z_1,a_2,\dots,z_2,\dots,a_k,\dots,z_k\le n'=n+1,$ $0\le k\le n'$, 
such that $$\{p_l,p_l',\dots,p_l+q_l\},\quad \{a_i,a'_i,\dots,z_i\}$$ have the disjoint/include property.
\end{cor}
\section{Evaluation}
For $\sigma\in S_n$, let $T^\sigma$ denote the Terada-$n$ 
$0\sigma(1)\dots\sigma(n)n'n''$. We would like to evaluate the sum
$$\dfrac{J_n}{n!}=\sum_{\sigma\in S_n}T\bullet T^\sigma;$$
when $T$ and $T^\sigma$ do not touch, their intersection number is $0$,
of course.
Proposition 1 and Corollary 1 imply that $(-)^nJ_n/n!$ is the sum of 
all the possible product
\begin{eqnarray}\prod_i(-)^{q_i}g^{q_i(q_i+1)/2}[p_ip'_i\cdots p_i+q_i]\cdot \prod_j[a_ja'_j\cdots z_j],\label{eq:3.1}\end{eqnarray}
where $1\le p_i,p_i+q_i\le n,\ 0\le a_j,\dots,z_j\le n+1,$ and $[a_ia'_i\cdots,z_i]\not=[01\cdots n']$ and,
$$\{p_i,p'_i,\dots,p_i+q_i\},\quad\{a_j,a'_j,\dots,z_j\}$$
have disjoint/include property. Empty products are regarded as $1$.
\par\medskip\noindent 
{\bf Notation:} For $1\le p,p+q\le n,$ put
$$\begin{array}{ll}
\langle pp'\cdots p+q\rangle&:=[pp'\cdots p+q](1+(-)^qg^{(q+1)q/2})\\[1mm]
&=\dfrac{1+(-)^qg^{q+1\choose2}}{g^{2{q+1\choose2}}-1}=\dfrac{-1}{1-(-)^qg^{q+1\choose2}}.\end{array}$$
\par\noindent 
{\bf Terminology:} $\langle\cdots\rangle$ and $[\cdots]$ are called 
{\it sequences}. A {\it monomial} in $0,1,\dots,n,n+1$ is a product of 
sequences
$$\prod_i\langle p_ip'_i\cdots p_i+q_i\rangle\cdot \prod_j[a_ja'_j\cdots z_j],$$ 
where $$\{p_i,p'_i,\dots,p_i+q_i\},\quad 1\le p_i,p_i+q_i\le n,$$
and
$$\{a_j,a'_j,\dots,z_j\},\quad 0\le a_j,\dots\le n+1,\quad 0=a_j\mbox{\ or\ } z_j=n+1$$
have disjoint/include property. The {\it length} of a monomial is the 
length of the longest sequence.
\par\smallskip
For  example, when $n=4$,
$$[0123]\langle23\rangle\langle123\rangle[45],\quad\langle12\rangle\langle34\rangle\langle1234\rangle[12345]$$
are monomilas of length 4 and 5.
\par\smallskip
Let $U$ be a family of subsets of $\{0,1,\dots,n'\}$ with disjoint/include property, and $J(U)$ the sum of possible products in the form (\ref{eq:3.1}), where$$U=\cup\{p_i,\dots,p_i+q_i\}\cup\cup\{a_j,\dots,z_j\}.$$
The sum $J(U)$ can be computed as follows: Divide $U$ into three subsets:
$$\begin{array}{ll}
{}_1U_n&:=\{\{p_i,\dots,p_i+q_i\}\in U\mid 1\le p_i,p_i+q_i\le n\},\\[1mm]
{}_0U&:=\{\{0,\dots,z_i\}\in U\mid z_i\le n\},\\[1mm]
U_{n'}&:=\{\{a_i,\dots,n'\}\in U\mid 1\le a_i\}.
\end{array}$$
Then the sum $J(U)$ can be factorized as
$$J(U)=\prod_{{}_1U_n}\langle p_i\cdots p_i+q_i\rangle\prod_{{}_0U}[0\cdots z_i]\prod_{U_{n'}}[a_i\dots n'].$$
We hope that the reader readily guess the reason of the above claim by the following example for $n=3$ and $U=\{\{1,2\},\{1,2,3\}\}$ (cf. \S4.2):
$$[12][123]+g[12]\cdot[123]-g^3[123]\cdot[12]+g[12](-g^3[123])$$
$$=(1-g)[12]\cdot(1+g^3)[123]
=\langle12\rangle\langle123\rangle.$$
Letting $U$ vary all the families of  subsets of $\{0,1,\dots,n'\}$ with disjoint/include property, we have
\begin{prp}
The quantity $(-)^nJ_n/n!$ is the sum of the possible monomials in 
$0,1,\dots,n,n+1$ of length at most $n+1$. 
\end{prp}
For each monomial, notice the longest sequence including $0$, and the longest sequence including $n+1$. Then we are led to
\begin{prp}For $k=1,\dots,n,$ put
$$\begin{array}{rl}
X_k&:=\mbox{sum of the possible monomilas in $1,2,\dots k$},\\
 A_k&:=[01\cdots k]\times\{\mbox{sum of the possible monomilas in $0,1,\dots k$ of length at most $k$}\},\\
B_k&:=[n-k+1\cdots n\  n+1]\times\{\mbox{sum of the possible monomilas} \\
&\quad\mbox{ in $n-k+1,\dots,n+1$ of length at most $k$}\}.\end{array}$$
Then we have
$$(-)^n\dfrac{J_n}{n!}=\sum_{i,j\ge0,\ i+j\le n}A_i\ X_{n-i-j}\ B_j.$$
\end{prp}
\subsection{Evaluation of $X_n$}
We found that the sum $X_n$ of the terms free of $a$ and $c$ is the sum of monomials in $1,\dots,n$.
Indeed we have
$$\begin{array}{rl}
X2&=1+\langle12\rangle,\\[1mm]
X3&=(1+\langle123\rangle)(1+\langle12\rangle+\langle23\rangle),\\[1mm]
X4&=(1+\langle1234\rangle)(1+\langle12\rangle+\langle23\rangle+\langle34\rangle+\langle12\rangle\langle34\rangle\\[1mm]
&+(1+\langle12\rangle+\langle23\rangle)\langle123\rangle+(1+\langle23\rangle+\langle34\rangle)\langle234\rangle),\\[1mm]
&\cdots
\end{array}$$
These can be evaluated as
$$\begin{array}{rl}
X_2&=\dfrac{g}{g+1},\\[1mm]
X_3&=\dfrac{g^3}{(g+1)(g^2+g+1)},\\[1mm]
X_4&=\dfrac{g^{6}}{(g+1)(g^2+g+1)(g^3+g^2+g+1)},\\[1mm]
&\cdots\end{array}$$
Let $Y(k,n)$ be the sum of the monomials in $1,\dots,n$ of length exactly $k$, and $X(k,n)$ be the sum of the monomials in $1,\dots,n$ of length at most $k$. Then we have
$$X(k,n)=\sum_{j=1}^kY(j,n),\quad X_n=X(n,n),$$ 
where we put $Y(1,2)=Y(1,3)=\cdots=1$, so we have $X(1,2)=X(1,3)=\cdots=1.$
\par\medskip\noindent
{\bf Notation:}
$$\begin{array}{ll}
[n]&:=[n]_g=1+g+g^2+\cdots+g^{n-1},\\[2mm]
{[n]!} &:=[n]_g!=[n][n-1]\cdots[1],\\[2mm]
\displaystyle{\left[\begin{array}{c}n\\m\end{array}\right]} &:=
\displaystyle{\left[\begin{array}{c}n\\m\end{array}\right]}_g=\dfrac{[n]!}{[m]!\ [m]!},\quad 0\le m\le n,\\[1mm]
(a)_n&:=(a;g)_n=(1-a)(1-ag)\cdots(1-ag^{n-1}).\phantom{\dfrac{A}{A}}
\end{array}$$
Since we have, by definition,
$$Y(n,n)=\langle1\cdots n\rangle\ X(n-1,n),\quad
X(n,n)=Y(n,n)+X(n-1,n)$$ and $$\langle1\cdots
n\rangle=\frac{-1}{1+(-)^ng^{n\choose2}},$$
\begin{lmm}
$$X(n-1,n)=\dfrac{g^{n\choose2}+(-)^n}{[n]!},\quad n=2,3,\dots$$\end{lmm}
implies
\begin{prp}
$$Y(n,n)=\frac{(-)^{n+1}}{[n]!},\quad X(n,n)=\frac{g^{n\choose2}}{[n]!}.$$\end{prp}
Proof. We prove this Lemma (as well as this Proposition) by induction on $n$.
For each monomial of length at most $n-1$ in $\{1,2,\dots,n\}$, we
notice the longest sequence including 1, and divide the sequences in two 
parts: those which are part of the longest one, and those which are
disjoint with the longest one. In this way, we have
$$\begin{array}{ll}
X(n-1,n)&=\langle1\dots n-1\rangle\ X(n-2,n-1)\\[1mm]
        &+\ \langle1\dots n-2\rangle\ X(n-3,n-2)\ X(2,2)\\
&\vdots\\
        &+\ \langle1\dots n-k\rangle\ X(k-1,k)\ X(n-k,n-k)\\
&\vdots\\
        &+\ \langle12\rangle\ X(12)\ X(n-2,n-2)+X(n-1,n-1).\end{array}
$$
By the induction hypothesis of the Lemma and the Proposition, we have
$$\langle1\dots n-k\rangle\ X(k-1,k)\ X(n-k,n-k)=\frac{(-)^{k+1}}{[k]!}\ 
\frac{g^{n-k\choose2}}{[n-k]!},$$
and so
$$X(n-1,n)=\sum_{k=1}^{n-1}\frac{(-)^{k+1}}{[k]!}\ 
\frac{g^{n-k\choose2}}{[n-k]!}=\frac1{[n]!}\sum_{k=1}^{n-1}(-)^{k+1}\displaystyle{\left[\begin{array}{c}n\\k\end{array}\right]}g^{n-k\choose2}.$$
On the other hand, the $g$-binomial theorem (cf. \cite{GR})
$$(x)_n=\sum_{k=0}^n\left[\begin{array}{c}n\\k\end{array}\right]g^{1+\cdots+(k-1)}(-x)^k$$
with $x=1$ yields
$$\sum_{k=0}^n\left[\begin{array}{c}n\\k\end{array}\right]g^{k\choose2}(-)^k=0\mbox{\quad
or\quad}\sum_{k=0}^n\left[\begin{array}{c}n\\k\end{array}\right]g^{n-k\choose2}(-)^{n-k}=0.$$
This leads to
$$X(n-1,n)=\frac{g^{n\choose2}+(-)^n}{[n]!},$$
proving the Lemma, and so the Proposition.
\subsection{Evaluation of $A_k$}
\begin{prp}$$\begin{array}{ll}
A_k&=[01\cdots k]\times\{\mbox{sum of the monomials in $0,1,\dots,k$ of length at most $k$}\}\\[1mm]
&=[01\cdots k]\displaystyle{\sum_{p=0}^{k-1}A_p\ X_{k-p}},\end{array}$$
where $A_0=1$.\end{prp}
For example
$$\begin{array}{ll}
A_1&=[01],\\[1mm]
A_2&=[012](X_2+[01]),\\[1mm]
A_3&=[0123](X_3+[01]X_2+A_2),\ \cdots\end{array}$$
\begin{prp}
$$A_k=\dfrac{(-)^k\ (1-g)^k}{(a)_k\ (g)_k},\quad k=0,1,\dots$$
$B_k$ is given from this formula by replacing $a$ by $b$.
\end{prp}
Proof. We compute
$$\sum_{k=0}^{n-1}A_kX_{n-k}=X_n\left[1+A_1\frac{X_{n-1}}{X_n}+\cdots+A_k\frac{X_{n-k}}{X_n}+\cdots+A_{n-1}\frac{X_{1}}{X_n}\right].$$
Since we have
$$\begin{array}{ll}
X_k=\dfrac{(1-g)^kg^{k\choose2}}{(g)_k}&=\dfrac1{(1+g^{-1})\cdots(1+g^{-1}+\cdots+g^{-k+1})}\\[1mm]
\phantom{\dfrac{\dfrac{A}{B}}{B}}&=\dfrac{(1-g^{-1})^k}{(1-g^{-1})(1-g^{-2})\cdots(1-g^{-k})}=\dfrac{(1-g^{-1})^k}{(g^{-k})_k},\end{array}$$
and so
$$\begin{array}{ll}
\dfrac{X_{n-k}}{X_n}&=(1-g^{-1})^{-k}\dfrac{(g^{-n})_n}{(g^{n-k})_{n-k}}\\[1mm]
\phantom{\dfrac{\dfrac{A}{B}}{B}}&=(1-g^{-1})^{-k}(1-g^{-n})\cdots(1-g^{-n+k-1})=\dfrac{(g^{-n})_k}{(1-g^{-1})^k}=g^k\dfrac{(g^{-n})_k}{(g-1)^k},
\end{array}$$
by the induction hypothesis, we have
$$A_k\frac{X_{n-k}}{X_n}=\frac{(g^{-n})_k}{(a)_k(g)_k}g^k.$$
On the other hand, the reversed $q$-Chu-Vandermonde formula (cf. \S6.2) on the basic hypergeometric 
function $\varphi$:
$$\sum_{k=0}^n\frac{(q^{-n})_k\ (b)_k}{(q)_k\ (c)_k}=:{}_2\varphi_1\left(\begin{array}{c}q^{-n},b\\c\end{array};q,q\right)=\dfrac{(c/b;q)_n}{(c;q)_n}\
b^n$$
with $q\to g,c\to a,$ and $b\to0$ yields
$$\sum_{k=0}^n\frac{(g^{-n})_k}{(a)_k(g)_k}\ g^k=\frac{(-a)^n\ g^{1+\cdots+(n+1)}}{(a)_n}=\frac{(-a)^n\ g^{n\choose2}}{(a)_n}.$$
Thus we have
$$\begin{array}{ll}
\displaystyle{1+A_1\frac{X_{n-1}}{X_n}+\cdots+A_k\frac{X_{n-k}}{X_n}+\cdots+A_{n-1}\frac{X_{1}}{X_n}}&=\dfrac{(-a)^n\ g^{n\choose2}}{(a)_n}-\dfrac{(g^{-n})_ng^n}{(a)_n(g)_n}\\
\phantom{\dfrac{\dfrac{A}{B}}{B}}&=\dfrac{a^ng^{n\choose2}-g^{-{n\choose2}}}{(a)_n}(-)^n,\end{array}$$
and so
$$\begin{array}{ll}
\sum_{k=0}^{n-1}A_kX_{n-k}&=X_n\left[1+\cdots+A_n\dfrac{X_1}{X_n}\right]\\[1mm]
\phantom{\dfrac{\dfrac{A}{B}}{B}}&=\dfrac{(1-g)^n\ g^{n\choose2}}{(g)_n}\cdot\dfrac{a^ng^{n\choose2}-g^{-{n\choose2}}}{(a)_n}(-)^n=\dfrac{(g-1)^n(a^ng^{2{n\choose2}}-1)}{(a)_n\ (g)_n}.\end{array}$$
Therefore we have
$$A_n=[01\cdots n]\sum_{k=0}^{n-1}A_kX_{n-k}=\dfrac{(g-1)^n}{(a)_n\ (g)_n},$$
ending the proof.

\subsection{Evaluation of $J_n$, end of the proof of Theorem 1}
\begin{prp}
$$\sum_{i,j\ge0,\ i+j\le n}A_i\ X_{n-i-j}\ B_j
=\dfrac{(g-1)^n\ (abg^{n-1})_n}{(a)_n\ (b)_n\ (g)_n}.
$$\end{prp}
Proof. In the course of the proof of the previous Proposition, we proved 
$$\sum_{i=0}^mA_iX_{m-i}=\frac{(g-1)^ma^mg^{2{m\choose2}}}{(a)_m\
(g)_m}=:C_m.$$
We compute
$$\sum_{j=0}^n\sum_{i=0}^{n-j}A_iX_{n-i-j}B_j=\sum_{j=0}^nC_{n-j}B_j=C_n\sum_{j=0}^nB_j\frac{C_{n-j}}{C_n}.$$
Since we have
$$\begin{array}{ll}(a)_n&=(1-a)(1-ag)\cdots(1-ag^{n-1})\\[1mm]
&=(-a)^ng^{1+\cdots+(n-1)}(1-a^{-1})(1-a^{-1}g^{-1})\cdots(1-a^{-1}g^{1-n}),\end{array}$$
and so
$$\frac{(a)_n}{(a)_{n-j}}=(-a)^j\frac{g^{1+\cdots+(n-1)}}{g^{1+\cdots+(n-j-1)}}(a^{-1}g^{1-n})_j,\quad
\frac{(g)_n}{(g)_{n-j}}=(-g)^j\frac{g^{1+\cdots+(n-1)}}{g^{1+\cdots+(n-j-1)}}(g^{-n})_j,$$
we have
$$\begin{array}{ll}
\dfrac{C_{n-j}}{C_n}&=(g-1)^{-j}a^{-j}\dfrac{(a)_n\ (g)_n}{(a)_{n-j}(g)_{n-j}}g^{(n-j)(n-j-1)g^{-n(n-1)}}\\
\phantom{\dfrac{A}{B}}
&=(a^{-1}g^{1-n})_j\ (g^{-n})_j\ (g-1)^{-j}g^j,\end{array}$$
and so
$$B_j\ \dfrac{C_{n-j}}{C_n}=\dfrac{(g^{-n})_j\ (a^{-1}g^{1-n})_j}{(b)_j\
(g)_j}g^j. $$
On the other hand the reversed $q$-Chu-Vandermonde formula above yields
$$\sum_{j=0}^n\dfrac{(g^{-n})_j\ (a^{-1}g^{1-n})_j}{(b)_j\ (g)_j}g^j
=\dfrac{(abg^{n-1})_n}{(b)_n}(a^{-1}g^{1-n})^n.$$
Now we finish the proof:
$$\begin{array}{ll}
\displaystyle{\sum_{j=0}^n\sum_{i=0}^{n-j}A_iX_{n-i-j}B_j}&=\displaystyle{C_n\sum_{j=0}^nB_j\frac{C_{n-j}}{C_n}}
=C_n\dfrac{(a{b}g^{n-1})_n}{(b)_n}(a^{-1}g^{1-n})^n\\ \phantom{\frac{\dfrac{A}{A}}{\dfrac{B}{B}}}
&=\dfrac{(g-1)^na^ng^{n(n-1)}}{(a)_n\ (g)_n}\cdot\dfrac{(abg^{n-1})_n}{(b)_n}(a^{-1}g^{1-n})^n\\ \phantom{\frac{\dfrac{A}{A}}{\dfrac{B}{B}}}
&=\dfrac{(g-1)^n\ (abg^{n-1})_n}{(a)_n\ (b)_n\ (g)_n}.\phantom{\frac{\dfrac{A}{A}}{\dfrac{B}{B}}}\end{array}$$

\section{Review of 2D and 3D cases}
\subsection{2D case}
The Pentagon $T:=xy=0xy1$ is a Terada-2 with code $01234$. It has \par\noindent
1D faces 
$$(01),\ (12),\ (23),\ (34)=(012),\ (40)=(123),$$
0D faces
$$(01)\cap(23),\ (01)\cap(34),\ (12)\cap(34),\ (12)\cap(40),\ (23)\cap(40).$$
\par\medskip\noindent
Exponents of the 1D faces are
$$(01)\ \cdots\  a,\quad (12)\  \cdots\  g^2,\quad (23)\ \cdots\  b,$$
$$(34)=(012)\ \cdots\  a^2g^{2}, \quad(40)=(123)\ \cdots g^{2}b^2;$$
put
$$[01]=\dfrac1{a-1},\quad [12]=\dfrac1{g^2-1},\quad [34]=[012]=\dfrac1{a^2g^2-1},$$
and so on. Then we have the self-intersection number of $T=xy=01234$:
$$\begin{array}{ll}xy\bullet xy&=1+[01]+[12]+[23]+[34]+[40]\\[2mm]
&+[01][23]+[01][34]+[12][34]+[12][40]+[23][40].\end{array}$$
Since $yx=02134$ is adjacent to $T$ along the edge $(12)$, computing the self-intersection number of the edge $(12)$, we have
$$xy\bullet yx=-g[12]\{1+[34]+[40]\}.$$
Let us add these:
$$\begin{array}{ll}1+(1-g)[12]&+\ [01]+[012](1+(1-g)[12]+[01])\\[2mm]
&+\ [23]+[123](1+(1-g)[12]+[23])+[01][23].\end{array}$$
Note that the sum of the terms free of $a$ and $b$ is $1+(1-g)[12]$, and 
the sum of the terms with the factor $a^2g^2-1$ is $[012](1+(1-g)[12]+[01])$.
\subsection{3D case}
The Terada-3 $T:=0xyz1$ is now coded by $012345$. It has \par\noindent
2D faces of 2 kinds
$$(01),\ \dots,\ (50);\quad (012)=(345),\ (123)=(456),\ (234)=(501),$$
1D faces of 2 kinds
$$(01)\cap(23);\quad (01)\cap(012),$$
0D faces of 2 kinds
$$(01)\cap(23)\cap(45),\ (12)\cap(34)\cap(50);\quad (01)\cap(12)\cap(012).$$
\par\medskip\noindent
Exponents of the 2D faces are
$$(01)\ \cdots\  a,\quad (12),(23)\ \cdots\  g^2,\quad (34)\ \cdots\  b,$$
$$(45)=(0123)\ \cdots\ a^3g^{6}, \quad(50)=(1234)\ \cdots\ g^{6}b^3,$$
$$(012)\ \cdots\ a^2g^2,\quad (123)\ \cdots\ g^6,\quad (234)\ \cdots\ g^2b^2.$$
Recall the notation:
$$[01]:=F(a),\ [12]:=F(g^2),\ \dots, \ [50]:=F(g^6b^3), \ \dots, \ [234]:=F(g^2b^2).$$
The faces touching {\it the} $3!$ T3's are $xyz=123$ itself and\par\noindent
through 2D faces: $$(12),\  (23);\quad (123),$$
through 1D faces: $$(12)\cap(123),\ (23)\cap(123).$$
Then counting all the faces of $012345$, we have the self-intersection number of $T$:
$$xyz\bullet xyz=-\{1+[01]+\cdots+[50]+[012]+[123]+[234]+ [01][23]+\cdots\}. $$
Through the face $(12)$, $T$ is adjacent to $yxz=021345$; evaluating the self-intersection number of the face $(01)$, we have
$$xyz\bullet yxz=g[12]\{1+[34]+[45]+[50]+[123]+[012]+[34][50]+[34][012]+\cdots\}$$
Through the face $(123)$, $T$ is adjacent to $zyx=032145$; evaluating the self-intersection number of the face $(123)$, we have
$$xyz\bullet zyx=-g^3[123]\{1+[12]+[23]+[45]+[50]+[12][45]+\cdots\}$$
Along the face $(12)\cap(123)$, $T$ touches $zxy=031245$; evaluating the self-intersection number of the face $(12)\cap(123)$, we have
$$xyz\bullet zxy=g[12]g^3[123]\{1+[45]+[50]\}.$$
Let us add these terms.
 First the sum $X_3$ of the terms free of $a$ and $b$: 
$$\begin{array}{ll}
X_3&:=1+[12]+[23]+[123]+[12][123]+[23][123]\\[2mm]
&\quad-g[12]\{1+[123]\}-g[23]\{1+[123]\}\ \ +g^3[123]\{1+[12]+[23]\}\\[2mm]
&\quad-g[12]g^3[123]\{1\}-g[23]g^3[123]\{1\}\\[2mm]
&=(1-[12](g-1)-[23](g-1))(1+[123](1+g^3))\\[2mm]
&=\dfrac{g^3}{(g+1)(g^2+g+1)}.\end{array}$$
Next sum of the terms with the factor $[0123]=(a^3g^6-1)^{-1}$ is $[0123]$ times
$$\begin{array}{l}\quad -X_3-([01]+[01][23]+[012](1+[01]+[12])+g[12][012]+g[23][01]\\[2mm]
=-X_3-[01](1-[23](g-1))-[012](1-[12](g-1)+[01])\\[2mm]
=\displaystyle{-X_3-\frac{1}{a-1}\frac{g}{g+1}-\frac1{a^2g^2-1}\frac{ag+1}{(a-1)(g+1)}}\\[2mm]
=\dfrac{-(a^3g^6-1)}{(a-1)(ag-1)(ag^2-1)(g+1)(g^2+g+1)}.\end{array}$$

\section{Proof of Theorem 2} Evaluation of intersection of cohomology is simpler than that of homology. Especially the evaluation of the self-intersection number of the form corresponding to a chamber is the sum of the contribution of every vertices. Let us see for example the $n$-beta function
$$B(\alpha_0,\dots,\alpha_n):=\int_{t_j>0,\ t_1+\cdots+t_n<1}
t_1^{\alpha_1-1}\cdots t_n^{\alpha_n-1}(1-t_1-\cdots-t_n)^{\alpha_0-1}
dt_1\dots dt_n.$$
The quadratic relation reads
$$B(\alpha_0,\dots,\alpha_n)\cdot B(-\alpha_0,\dots,-\alpha_n)=
\displaystyle{\frac{1-\prod_0^na_i}{\prod_0^n(1-a_i)}\cdot(2\pi i)^n\frac{\alpha_0+\cdots+\alpha_n}{\alpha_0\cdots\alpha_n}},$$
where the first factor of the right hand-side is the self-intersection of the loaded simplex
$$\Delta:\quad t_j>0,\ t_1+\cdots+t_n<1$$
and the second factor is the self-intersection of the twisted form
$$\frac{dt_1\dots dt_n}{t_1\cdots t_n(1- t_1-\cdots-t_n)}.$$
The former is the sum of the contributions of {\it all faces} of $\Delta$ :
$$1+\sum_i\frac1{a_i-1}+\sum_{i<j}\frac1{(a_i-1)(a_j-1)}+\cdots+\sum_i\frac1{\prod_{j\not=i}(a_j-1)},$$
and the latter is the sum of the contribution of just the {\it vertices} of $\Delta$:
$$\sum_i\frac1{\prod_{j\not=i}\alpha_j}\quad\mbox{at the point $x_j=0 \ (j\not=i)$}.$$
\par\medskip
Now we consider Terada-$n$ $T$ as a polyhedron. The self-intersection number
$$\omega:=\frac{dt_1\wedge\cdots\wedge dt_n}{\prod_{i=1}^n\ t_i\ (1-t_i)}$$
is $(2\pi i)^n\ n!$ times the sum of the local contribution at every vertex of $T$ not touching another Terada-$n$, because $\omega$ has no poles of type
$$\frac{dt_1}{t_1}\wedge\cdots\wedge\frac{dt_n}{t_n}$$
at vertices touching another Terada-$n$. Such a vertex will be said to be 
admissible. Let us work in this line. 
\begin{itemize}
\item 2D case: admissible vertices are (in Figure 1, they are marked)
$$(012)\cap(01),\  (01)\cap(34),\  (34)\cap(234);$$ so the sum of the contributions are 
$$\frac1{2\alpha+2\g}\ \frac1{\alpha}+\frac1\alpha\ \frac1\beta+\frac1\beta\ \frac1{2\beta+2\g}=\frac12\ \frac{(\alpha+\beta+\g)(\alpha+\beta+2\g)}{\alpha(\alpha+\g)\beta(\beta+\g)}.$$
\item 3D case: admissible vertices are (in Figure 3, they are marked)
$$(0123)\cap(012)\cap(01),\ (012)\cap(01)\cap(45),\ (01)\cap(45)\cap(345),\ (45)\cap(345)\cap(2345);$$ so the sum of the contributions are 
$$\frac1{3\alpha+6\g}\ \frac1{2\alpha+2\g}\ \frac1\alpha+ 
                       \frac1{2\alpha+2\g}\ \frac1\alpha\ \frac1\beta+ 
                                            \frac1\alpha\ \frac1\beta\ \frac1{2b+2g}+ 
                                                      \frac1\beta\ \frac1{2\beta+2\g}\ \frac1{3\beta+6\g}$$$$=\frac1{3!}\frac{(\alpha+\beta+2\g)(\alpha+\beta+3\g)(\alpha+\beta+4\g)}{\alpha(\alpha+\g)(\alpha+2\g)\beta(\beta+\g)(\beta+2\g)}.
$$
\item $n$D case: admissible vertices are, for $k=0,\dots,n,$
$$(01)\cap(012)\cap\cdots\cap(01\cdots k)\cap(k'\cdots n')\cap\cdots\cap({}'nnn')\cap(nn');$$
so the sum of the contributions are
$$\sum_{k=0}^n\frac1{\prod_{i=1}^k(i\alpha+i(i-1)\g)\ 
\prod_{j=1}^{n-k}(j\beta+j(j-1)\g)}.$$
\end{itemize}
Theorem 2 claims that this is equal to
$$=\frac1{n!}\frac{(\alpha+\beta+(n-1)\g)\ (\alpha+\beta+n\g)\cdots(\alpha+\beta+(2n-2)\g)}{\alpha(\alpha+\g)\cdots(\alpha+(n-1)\g)\ \beta(\beta+\g)\cdots(\beta+(n-1)\g)}.$$ 
Proof. In this proof $(\alpha)_n:=\alpha(\alpha+1)\cdots(\alpha+n-1).$ The sum in question equals
$$\begin{array}{ll}
&=\displaystyle{\dfrac1{n!}\sum_{k=0}^n{n\choose k}\dfrac1{\prod_{i=1}^k(\alpha+(i-1)\g)\prod_{j=1}^{n-k}(\beta+(j-1)\g)}}\\
&=\displaystyle{\dfrac1{n!}\dfrac1{\prod_{j=1}^n(\beta+(j-1)\g)}\sum_{k=0}^n{n\choose k}\dfrac{\prod_{j=n-k+1}^n(\beta+(j-1)\g)}{\prod_{i=1}^k(\alpha+(i-1)\g)}}\\
&=\displaystyle{\dfrac1{n!}\dfrac1{\prod_{j=1}^n(\beta+(j-1)\g)}
\sum_{k=0}^n\dfrac{(-n)_k\left(-\dfrac{\beta}{\g}-n+1\right)_k}{k!\ \left(\dfrac{\alpha}{\g}\right)_k}}.
\end{array}$$
On the other hand, the Chu-Vandermonde formula (a special case of Gauss's summation formula) for the value of the hypergeometric function ${}_2F_1$ at $1$:
$${}_2F_1\left(\begin{array}{c}-n,\beta\\ \g\end{array};1\right)=\dfrac{(\gamma-\beta)_n}{(\g)_n}$$
implies that the sum above equals 
$$\sum_{k=0}^n\dfrac{(-n)_k\left(-\dfrac{\beta}{\g}-n+1\right)_k}{k!\ \left(\dfrac{\alpha}{\g}\right)_k}
=\dfrac{\prod_{j=1}^n(\alpha+\beta+(n-1+j-1)\g)}{\prod_{j=1}^n(\alpha+(j-1)\g)},$$
which ends the proof.

\section{Appendix}
\subsection{Standard loading} For a set of hyperplanes $f_j(t)=0$ in the complex affine $n$-space and complex numbers (called exponents) $\alpha_j$, we consider a (multi-valued) function
$$u=\prod_jf_j(t)^{\alpha_j}$$
on $X:=\C^n-\cup_j\{f_j=0\}.$ The function $u$ determines the local system 
$\mathcal S=\mathcal S_u$ on $X$.  When every $f_j$ is defined over $\R$, the real locus $X_\R$ breaks into simply connected chambers. For each chamber $D$ we load 
$$\prod_j(\varepsilon_j\cdot f_j(t))^{\alpha_j},\qquad \arg\varepsilon_j\cdot f_j(t)=0, $$
where $\varepsilon_j=\pm$ is so determined that $\varepsilon_j\cdot f_j(t)$ is positive on $D$. This way of loading is said to be {\it standard}.
This loading is already used in \cite{MiY}, though the terminology `standard' is not used.

In all the other publications of the authors, chambers are loaded in a different way: each chamber is loaded with a branch (result of analytic continuation along a fixed path) of a fixed function element of $u$ at a specific base point.  This loading is often convenient in practical computation. But it has a fatal disadvantage: it depends on the choice of the base point and the paths. As a result, the intersection matrices are not symmetric.

On each chamber, the old loading and the standard loading differ only by a multiplicative constant, and their global monodromy representation (of the fundamental group) coincides, of course. So the differnce might look minor, but this standard loading has apparent advantages: it does not require any base point nor fixed paths, and the intersection matrices are symmetric.
\par\smallskip
Here we present the simplest example. In the $x$-space, we consider
$$u=x^\alpha(1-x)^\beta(2-x)^\delta.$$
Let us load the interval $I_0:=(0,1)$ and $I_1:=(1,2)$ in a standard way, that is, we load $I_0$ with $u_0:=x^\alpha(1-x)^\beta(2-x)^\delta$, and $I_1$ with $u_1:=x^\alpha(x-1)^\beta(2-x)^\delta$, where the aregument of positive numbers are supposed to be 0. By making regularizations (by putting kintamas), we evaluate the intersection numbers as:
$$({\rm reg\ } I_1)\bullet \check I_{0}=({\rm reg\ } I_{0})\bullet \check I_1=\frac{e^{\pi i\beta}}{e^{2\pi i\beta-1}}.$$
If we choose a base point in the lower half $x$-space, and load these intervals with the result of analytic continuations of u along paths in the lower half-space connecting the base point and the intervals, and if we write the loads on the intervals $I_0$ and $I_1$ by $u'_0$ and $u'_1$, respectively, they are related as
$$u'_0=Cu_0,\quad u'_1=Ce^{\pi i\beta} u_1,$$
where $C$ is a constant depending on the choice of the branch of $u$ at the base point. 
Their intersection numbers are now evaluated as 
$$({\rm reg\ } I_1)\bullet \check I_{0}=\frac{e^{2\pi i\beta}}{e^{2\pi i\beta}-1},\quad
({\rm reg\ } I_{0})\bullet \check I_1=\frac{{1}}{e^{2\pi i\beta}-1}.$$
These non-symmetric formulae are commonly presented, for instance in \cite{Yo1} (Chap 4 \S7).
\subsection{$q$-Chu-Vandermonde formulae}
We give a proof of the reversed $q$-Chu-Vandermonde formula
\begin{eqnarray}{}_2\varphi_1\left(\begin{array}{c}q^{-n},b\\c\end{array};q,q\right)=\dfrac{(c/b)_n}{(c)_n}\
b^n.\label{eq:invqCV}\end{eqnarray}
Reversing the order of the finite sum
$${}_2\varphi_1\left(\begin{array}{c}q^{-n},b\\c\end{array};q,x\right)=\sum_{i=0}^n A_i,\quad A_i=\dfrac{(q^{-n})_i(b)_i}{(q)_i(c)_i}x^n,$$
as
$$A_n\left(1+\dfrac{A_{i-1}}{A_n}+\cdots+\dfrac{A_i}{A_n}+\cdots+\dfrac{A_n}{A_n}\right),$$
we have
\begin{eqnarray}{}_2\varphi_1\left(\begin{array}{c}q^{-n},b\\c\end{array};q,x\right)=(-)^nq^{-{n+1\choose k}}\dfrac{(b)_n}{(c)_n}x^n\ {}_2\varphi_1\left(\begin{array}{c}q^{-n},c^{-1}q^{1-n}\\b^{-1}q^{1-n}\end{array};q,\dfrac{c}{b}\dfrac{q^{n+1}}{x}\right).\label{eq:inv}\end{eqnarray}
Indeed since 
$$\dfrac{(q^{-n})_{n-i}}{(q^{-n})_n}\cdot\dfrac{(q)_{n}}{(q)_{n-i}}=q^{(n+1)i}\dfrac{(q^{-n})_{i}}{(q)_i}\mbox{\quad and\quad } \dfrac{(c)_n}{(c)_{n-i}}\cdot\dfrac{(b)_{n-i}}{(b)_n}=\left(\dfrac{c}{b}\right)^i\frac{(c^{-1}q^{1-n})_i}{(b^{-1}q^{1-n})_i}$$
hold, we have
$$\frac{A_{n-i}}{A_n}=\frac{(q^{-n})_{n-i}(q)_n}{(q)_{n-i}(q^{-n})_n}\cdot\frac{(q)_{n-i}(c)_n}{(c)_{n-i}(b)_n}\ x^{-i}=\frac{(c^{-1}q^{1-n})_i}{(b^{-1}q^{1-n})_i}\dfrac{(q^{-n})_{i}}{(q)_i}\left(\dfrac{c}{b}\dfrac{q^{n+1}}{x}\right)^i,$$
which proves the above formula (\ref{eq:inv}). 
Putting $x=cq^n/b$ in this formula, we have
\begin{eqnarray}{}_2\varphi_1\left(\begin{array}{c}q^{-n},b\\c\end{array};q,\frac{cq^n}{b}\right)=(-)^n\left(\dfrac{c}{b}\right)^nq^{\frac12n(n-1)}\frac{(b)_n}{(c)_n}\ {}_2\varphi_1\left(\begin{array}{c}q^{-n},c^{-1}q^{1-n}\\b^{-1}q^{1-n}\end{array};q,q\right).\label{eq:inv2}\end{eqnarray}
On the other hand the $q$-analogue of the Gauss's sum due to Heine (cf. \cite{GR}) reads
$${}_2\varphi_1\left(\begin{array}{c}a,b\\c\end{array};q,\dfrac{c}{ab}\right)=\dfrac{(c/b)_\infty(c/b)_\infty}{(c)_\infty(c/ab)_\infty}.$$
This yields ($a=q^{-n}$) the $q$-Chu-Vandermonde formula
\begin{eqnarray}{}_2\varphi_1\left(\begin{array}{c}q^{-n},b\\c\end{array};q,\dfrac{c}{b}q^n\right)=\dfrac{(c/b)_n}{(c)_n}.\label{eq:qF1}\end{eqnarray}
Since the left hand-sides of (\ref{eq:inv2}) and (\ref{eq:qF1}) coincide, we have
$${}_2\varphi_1\left(\begin{array}{c}q^{-n},c^{-1}q^{1-n}\\b^{-1}q^{1-n}\end{array};q,q\right)=\dfrac{(c/b)_n}{(b)_n}(-)^n\left(\frac{b}{c}\right)^nq^{-\frac12n(n-1)}.$$
If we put $c^{-1}q^{1-n}\to b$ and $b^{-1}q^{1-n}\to c$, then this leads to (\ref{eq:invqCV}).
\par\bigskip\noindent
{\it Acknowledgement.} The second author is grateful to Joichi Kaneko for his encouragement.

\end{document}